\documentclass[12pt]{article}
 \usepackage{amsfonts}
\usepackage{amsmath}
\usepackage{amssymb}
\usepackage{graphicx}
\usepackage[usenames]{color}

\def\one{\mbox{$1\!\!\!\,\rule{0,1mm}{2,4mm}\,$}}

\def\mk{\medskip }

\def\cl{\mbox{\rm cl\,}}
\def\bit{\begin{itemize}}
\def\eit{\end{itemize}}
\def\vsp{\vspace*{1,5mm}\\ }
\def\n{\noindent }
\def\bk{\bigskip }

\def\E{\mathbb{E}}
\def\rr{\mathbb{R}}
\def\PP{\mathbb{P}}
\def\nn{\mathbb{N}}

\def\ct{continuous}
\def\oo{{\omega}}
\def\ooo{{\Omega}}
\def\<{\left<}
\def\>{\right>}
\def\({\left(}
\def\){\right)}
\def\ff{\forall }
\def\D{\Delta }
\def\9{{\infty}}
\def\barr{\begin{array}}
\def\earr{\end{array}}
\def\ov{\overline}
\def\wt{\widetilde}
\def\vp{{\varepsilon}}

\def\pp{{\partial}}
\def\dd{\displaystyle}
\def\a{{\alpha}}
\def\b{{\beta}}
\def\vf{{\varphi}}
\def\lbb{{\lambda}}
\def\g{{\gamma}}
\def\pas{\mathbb{P}\mbox{-a.s.}}
\def\de{{\delta}}

\def\calf{{\mathcal{F}}}
\def\calh{\mathcal{H}}
\def\calo{\mathcal{O}}
\def\calf{\mathcal{F}}

\def\vspp{\vspace*{1,5mm}\\ }
\def\n{\noindent }
\def\na{{\nabla}}

\def\inr{\mbox{ in }}
\def\onr{\mbox{ on }}
\def\ifr{\mbox{ if }}
\def\fwg{following}
\def\eq{equation}

\def\beq#1{\begin{equation}\label{e#1}}
\def\eeq{\end{equation}}
\def\1{^{-1}}

\newtheorem{theorem}{Theorem}[section]
\newtheorem{corollary}[theorem]{Corollary}
\newtheorem{definition}[theorem]{Definition}

\newtheorem{lemma}[theorem]{Lemma}

\newtheorem{proposition}[theorem]{Proposition}
\newtheorem{remark}[theorem]{Remark}

\title{{\bf Stochastic variational inequalities and~applications to~the~total variation flow perturbed by~linear multiplicative noise}}

\author{{\bf Viorel Barbu}\thanks{Octav Mayer Institute of Mathematics (Romanian Academy), 700506 Ia\c si, Romania (vb41@uaic.ro). This author was partially supported by a grant of the Romanian National Authority for Scientific Research,  CNCS-UEFISCDI Project PN-II-ID-PCE-2011-3-0027 and the BiBoS -- Research Center.}
\and {\bf Michael R\"ockner}\thanks{Fakult\"at f\"ur Mathematik, Universit\"at Bielefeld, D-33501 Bielefeld, Germany (roeckner@math.uni-bielefeld.de). This   research was also supported by the DFG through CRC 701.}}
\date{}

\begin{document}
\maketitle

\begin{abstract}
In this work, we introduce a new method to  prove the existence and uniqueness of a variational solution to the stochastic nonlinear diffusion equation
 $dX(t)={\rm div}\ \left[\frac{\nabla X(t)}{|\nabla X(t)|}\right]dt+X(t)dW(t)\mbox{ in }(0,\9)\times\calo,$ where $\calo$ is a bounded and open domain in $\rr^N$, $N\ge 1$, and $W(t)$ is a
Wiener process of the form $W(t)=\sum^\9_{k=1}\mu_ke_k\b_k(t)$, $e_k\in C^2(\ov\calo)\cap H^1_0(\calo),$ and $\b_k$, $k\in\nn$, are independent Brownian motions.  This is a sto\-chas\-tic diffusion equation with a highly
singular diffusivity term and one  main result established here is that, for all initial conditions in $L^2(\calo)$, it is well posed in a class of continuous  solutions to the corresponding stochastic  variational inequality. Thus one obtains a stochastic version of the (minimal) total variation flow. The new approach developed  here also allows to prove the finite
time extinction  of solutions in dimensions \mbox{$1\le N\le3$,} which is another main result of this work.\\ 
\n{\bf Keywords:} stochastic diffusion equation, Brownian motion, bounded variation, convex functions, bounded variation flow.
\end{abstract}

\section{Introduction}

We are concerned here with the  stochastic nonlinear diffusion  equation
\begin{equation}\label{e1.1}
\barr{l}
dX(t)={\rm div}[{\rm sgn}(\na X(t))]dt+X(t)dW(t)\quad\inr\ (0,\9)\times\calo,\vspp
X=0\quad\onr\ (0,\9)\times\pp\calo,\vspp
X(0)=x\quad\inr\ \calo,\earr \end{equation} where $\calo$ is a bounded and convex open domain in $\rr^N$, $N\ge 1$, with smooth boun\-dary $\pp\calo$ and $W(t)$ is a Wiener process of the form
\begin{equation}\label{e1.2}
W(t)=\sum^\9_{k=1}\mu_ke_k\b_k(t),\ \ t\ge0,\ \inr\ \calo,
\end{equation}
where $\mu_k$ are real numbers, $e_k\in C^2(\ov\calo)\cap H^1_0(\calo)$ forming an orthonormal basis in $L^2(\calo)$ and $\{\b_k\}^\9_{k=1}$ are independent Brownian motions on a stochastic basis $\{\ooo,\calf,\calf_t,\PP\}$. For simplicity, let us assume that $e_k$, $k\in\nn$,  are the eigen\-func\-tions of the Dirichlet Laplacian:
$$\mbox{$-\D e_k=\lbb_ke_k$ in $\calo;$ $e_k=0$ on $\pp\calo,$}$$
(but cf. Remark \ref{r11} (iii) below).

Throughout the paper, we assume
\bit\item[(H1)] $C^2_\9:=\dd\sum^\9_{k=1}\mu^2_k|e_k|^2_\9<\9,$ \eit
 and
\bit\item[(H2)] $D_\9:=\dd\sum^\9_{k=1}\mu_k|\na e_k|_\9<\9,$\eit
where $|\cdot|_\9$ denotes supremum norm in  $C(\ov\calo)$.

Define \begin{equation}
\label{12prim}
\mu(\xi):=\dd\sum^\9_{k=1}\mu^2_k e^2_k(\xi),\ \xi\in\calo.\end{equation}



The multi--valued graph ${\rm sgn}:\rr^N\to 2^{\rr^N}$ is defined by
\begin{equation}\label{e1.3}
{\rm sgn}\ r=r|r|^{-1}\ \ifr\ r\ne0;\ \ {\rm sgn}\ 0=\{r\in\rr^N;\ |r|\le1\},
\end{equation}
and $|\cdot|$ is the Euclidean norm of $\rr^N$. By the same symbol $|x |$, we shall denote the absolute value of $ x\in \rr$.
It should be emphasized that the homogeneous boundary condition arising in \eqref{e1.1} is in a certain sense formal because
 \eqref{e1.1} is not well posed in the classical Sobolev spaces with zero trace on the boun\-dary.

  In nonlinear diffusion theory,  equation \eqref{e1.1} is derived from the continuity equation perturbed by a Gaussian process proportional to the density $X(t)$ of the material, that is,
$$dX(t)=J(\na X(t))dt+X(t)dW(t),$$
where $J={\rm sgn}$ is the flux of the diffusing material.  (See \cite{12a}, \cite{22a}, \cite{14a}.)

Equation \eqref{e1.1} is   also  relevant as a mathematical model for  faceted  crystal growth under a stochastic perturbation as well as in material sciences (see \cite{10} for the deterministic model and
complete references on the  subject).   As a matter of fact, these models are based on differential gradient systems corresponding to a convex and nondifferentiable potential (energy).

Other recent applications refer to the PDE approach to image recovery (see, e.g., \cite{8} and also \cite{5}, \cite{7}). In fact, if $x\in L^2(\calo)$ is the blurred image, one  might  find  the restored image via the  total variation  flow $X=X(t)$ generated by the stochastic equation
\begin{equation}\label{e1.4}
\barr{l}
dX(t)={\rm div}\(\dd\frac{\na X(t)}{|\na X(t)|}\)dt+X(t)dW(t)\quad\inr\ (0,\9)\times\calo,\vspp
X(0)=x\ \inr \calo.\earr
\end{equation}In its deterministic form, this is the so-called {\it total variation based image restoration model} and its stochastic version \eqref{e1.4} arises naturally in this context  as perturbation of the {\it total variation flow} by a Gaussian (Wiener) noise (which explains the title of the paper).

It should be said that, due to its high singularity, equation \eqref{e1.1} does not have a solution in the standard sense for every initial condition in $L^2(\calo)$, that is, as an It\^o integral equation,  and this happens in the deterministic case, too.  However, this equation has a natural formulation in the framework of stochastic variational inequalities (SVI) (see Definition \ref{d3.1} below) and, as we show later on, it is well posed in this generalized sense. Below, we shall call solutions to such (SVI) {\it variational solutions} and solutions to standard It\^o-integral equations, as e.g. the solutions to the approximating equation \eqref{e1.6} below (see Proposition \ref{l4.2} (i), {\it ordinary variational solutions}.

 In \cite{3},   a complete existence and uniqueness result was proved for variational solutions to \eqref{e1.1} in the case of additive noise, that is,
 \begin{equation}\label{e1.5}
\barr{l}
dX(t)-{\rm div}[{\rm sgn}(\na X(t))]dt=dW(t)  \quad\inr\ (0,\9)\times\calo,\vspp
X(0)=x\ \inr \calo,\qquad X(t)=0\ \onr\ (0,\9)\times\pp\calo,\earr
\end{equation}
  if $1\le N\le2$.    For the multiplicative noise $X(t)dW(t)$, only the exis\-tence of a variational solution was proved and
uniqueness remained open. (See, however, the work \cite{10a} for recent results on this line, if $x\in H^1_0(\calo)$.)

In this paper, we prove the existence and uniqueness of variational solutions to \eqref{e1.1} in all dimensions $N\ge1$ (see Theorem \ref{t3.1}) and  all initial conditions $x\in L^2(\calo)$. We would like to stress that one main difficulty is when $x\in L^2(\calo)\setminus H^1_0(\calo)$, while the case $x\in H^1_0(\calo)$ is more standard (see Remark \ref{r3.6} below). Furthermore, we   prove the finite-time extinction of solutions with positive probability, if $N\le3$.

The approach we use here to prove the existence and uniqueness of \eqref{e1.1} is obtained approximating equation \eqref{e1.1} by
\begin{equation}\label{e1.6}
\barr{l}
dX-\D\wt\psi_\lbb(X)dt=X\,dW\ \inr\ (0,T)\times\calo,\vsp
X=0\ \onr\ (0,T)\times\pp\calo,\
X(0)=x\ \inr\ \calo,\earr\end{equation}
where $\wt\psi_\lbb(r)=\psi_\lbb(r)+\lbb r$ and $\psi_\lbb$ is the Yosida approximation of the graph \eqref{e1.3}.
By the substitution  $Y=e^{-W}X$ (``scaling``), we reduce \eqref{e1.1}  and \eqref{e1.6} to a random nonlinear diffusion equation   (see  \eqref{e4.1} and cf. \cite{3aa}, \cite{4}, \cite{6a}) and, again, we reformulate this random equation as a (this time, deterministic) variational inequality (VI), but with random coefficients (see Definition \ref{d4.1}). This equivalent formulation of \eqref{e1.1} (respectively \eqref{e1.6}) as a random partial differential equation (PDE) is crucial for the uniqueness proof of variational solutions to \eqref{e1.1} (see Section 5) and allows to obtain sharper regularity results for \eqref{e1.6}  (see, e.g., Proposition \ref{l4.2}(iii) and Lemma \ref{l5.4nou}) than those   obtained by a direct analysis of the stochastic equation  as in \cite{3}, \cite{6aa}. This approach which combines the analysis of approximating stochastic equations in connection with their  equivalent random deterministic PDE versions is by our knowledge new in the ge\-ne\-ral theory of stochastic PDE and represents one principal contribution of~this~work.

\section{Preliminaries}
\setcounter{equation}{0}

For every $1\le p\le\9$, by $L^p(\calo)$ we denote the space of all Lebesgue $p$-integrable functions on $\calo$ with  norm  $|\cdot|_p$.
The scalar product in $L^2(\calo)$ is denoted by $\<\cdot,\cdot\>$. $W^{1,p}(\calo)$ denotes the standard Sobolev space
$\{u\in L^p(\calo);$ $ \na u\in L^p(\calo)\}$ with the corresponding norm
 $$\|u\|_{1,p}:=\(\int_\calo|\na u|^p d\xi\)^{1/p},$$where $d\xi$ denotes the Lebesgue measure on $\calo$. $W^{1,p}_0(\calo)$ denotes the space $\{u\in W^{1,p}(\calo);\ u=0$ on $\pp\calo\}$.  We set $H^1_0(\calo)=W^{1,2}_0(\calo),$ $\|\cdot\|_1=\|\cdot\|_{1,2}$ and $H^2(\calo)=\{u\in L^2(\calo):D^2_{ij}u\in L^2(\calo),\ 1\le i,j\le N\}$, with its usual norm $\|\cdot\|_{H^2(\calo)}.$ $H^{-1}(\calo)$  with     norm $\|\cdot\|_{-1}$  denotes
the dual of $H^1_0(\calo)=W^{1,2}_0(\calo)$. By $BV(\calo)$ we denote the space of functions $u$ of bounded variation on $\calo$ and by $\|Du\|$ the variation of $u$, that is,
\begin{equation}\label{e2.1}\|Du\|=\sup\left\{\int_\calo u\ {\rm div}\ \vf d\xi ;\ \vf\in C^\9_0(\calo;\rr^N),\ |\vf|_\9\le1\right\}.\end{equation}
By $BV^0(\calo)$ we denote the space of the functions $u\in BV(\calo)$ with vanishing trace on $\pp\calo.$

Consider the function $\phi_0:L^1(\calo)\to\ov \rr=(-\9,+\9]$
$$\phi_{0}(u)=\left\{\barr{ll}
\|Du\|&\ifr\ u\in BV^0(\calo),\vspp
+\9&\mbox{ otherwise},\earr\right.$$ and  denote by $\cl \phi_0$ the lower semicontinuous closure of $ \phi_{0}$ in $L^1(\cal O)$, that~is,
\begin{equation}\label{e2.3} \cl \phi_0(u)= \inf\left\{\liminf \phi_0(u_n);\ u_n\to u\in L^1(\calo)\right\}.\end{equation}
As in \cite[p.~437]{1a} define, for $u\in L^1(\calo)$,
$$G(u)=\left\{\barr{ll}
\dd\int_\calo|\na u|d\xi&\mbox{ if }u\in W^{1,1}_0(\ooo),\vsp
+\9&\mbox{ otherwise}.\earr\right.$$

\n Then (e.g., by \cite[Theorem 3.9]{1}) it is easy to see that
$$\cl \phi_0=\cl\,G.$$
Hence, by \cite[Proposition 11.3.2]{1a}, for $u\in L^1(\calo)$,
 $$\cl \phi_0(u)=\left\{
\barr{ll}
 \|Du\| + \dd\int_{\partial{\calo}}|\gamma_0(u)| d \calh^{N-1}&\mbox{ if } u\in BV(\calo),\vsp
   +\infty&\mbox{ otherwise,}  \earr\right.$$
where $\gamma_0(u)$ is the trace of $u$ on the boundary and $d \calh^{N-1} $ is the Hausdorff measure.

 Let $\phi$ denote the restriction of $\cl \phi_0(u)$ to $L^2(\calo)$, i.e.,
 \begin{equation}\label{e2.4}
\barr{ll}
\phi(u)=
 \|Du\| + \dd\int_{\partial{\calo}}|\gamma_0(u)| d \calh^{N-1}&\mbox{ if } u\in BV(\calo)\cap L^2(\calo),\vsp
  \phi(u)= +\infty&\mbox{ if }u\in L^2(\calo) \setminus BV(\calo).\earr\end{equation}

\n By $\pp\phi:D(\pp\phi)\subset L^2(\calo)\to L^2(\calo)$,  we denote the subdifferential of $\phi$, that~is,
\begin{equation}\label{e2.5}\pp\phi(u)=\{\eta\in L^2(\calo);\ \phi(u)-\phi(v)\le\<\eta,u-v\>,\ \ff v\in D(\phi)\},
\end{equation}where
$$D(\phi)=\{u\in L^2(\calo);\ \phi(u)<\9\}=BV(\calo)\cap L^2(\calo).$$ It turns out (see \cite{1}) that $\eta \in\pp\phi(u)$ iff there is $z\in L^\9(\calo;\rr^N)$ such that $\eta=-{\rm div}\ z,$ $|z|_\9\le1$,
and $\int_\calo \eta ud\xi=\phi(u)$.(Here and everywhere in the following the derivatives are taken in the sense of distributions on $\calo.$)
The mapping $\pp\phi$ is not everywhere defined on $D(\phi)$, but it is maximal monotone in $L^2(\calo)$ and so generates a semigroup flow $u(t,x)=e^{-t\pp\phi}x$ which is the solution to the evolution
equation (see \cite{6}, p.~72, \cite{2}, p.~47)
\begin{equation}\label{e2.6}\frac{du}{dt}\,(t)+\pp\phi(u(t))\ni0,\ \ \ff t\ge0,\ u(0)=x,
\end{equation}for each $x\in\ov{D(\phi)}=L^2(\calo).$ More precisely, for $x\in L^2(\calo)$,  there is a unique strong  solution $u:[0,\9]\to L^2(\calo)$  to \eqref{e2.6}   and, for each $T>0$,
\begin{eqnarray}
&&\sqrt{t}\ \frac{du}{dt}\in L^2(0,T;L^2(\calo)),\ t\phi(u(t))\in L^\9(0,T),\phi(u)\in L^1(0,T),\label{e2.7}\\[2mm]
&&t\,\frac{du}{dt}\in L^\9(0,T;L^2(\calo)),\ u\in C([0,T];L^2(\calo)). \label{e2.8}
\end{eqnarray}(See \cite{2}, p. 158.) In fact, if $u\in W^{1,1}_0(\calo)$ and $\eta\in {\rm div}\,[{\rm sgn}(\na u)]\cap L^2(\calo)\ne\emptyset$, then it is easily seen that $u\in D(\pp\phi)
$ and $ \eta\in \pp\phi(u).$

 We can rewrite equation \eqref{e1.1} as
\begin{equation}\label{e2.9}\barr{l}
dX(t)+\pp\phi(X(t))dt\ni X(t)dW(t),\ \ t\ge0,\vspp
X(0)=x.\earr
\end{equation}However, since the multi-valued mapping $\pp\phi:L^2(\calo)\to L^2(\calo)$ is highly singular, at present no general existence result for stochastic infinite dimensional
equations of subgradient type is applicable to the present situation and so a direct approach should be used in order to get existence and uniqueness of solutions for \eqref{e2.9}.

In the \fwg, $L^p(0,T;E)$, $1\le p\le\9$, and $E$ a Banach space, denotes the space of all
Bochner measurable functions $u:(0,T)\to E$ with $\|u\|_E\in L^p(0,T).$  By $C([0,T];E)$ we denote the space of all the continuous $E$-valued functions on $[0,T]$.  We also use
the notation
$$W^{1,p}([0,T];E)=\left\{u\in L^p(0,T;E),\ \dd\frac{du}{dt}\in L^p(0,T;E)\right\},$$ where $\dd\frac{du}{dt}$ is taken in sense of $E$-valued distributions on $(0,T)$.
(We recall that any $u\in W^{1,p}([0,T];E)$ is absolutely continuous and a.e. differentiable.)

  The plan of the rest of the paper is the following.

In Section 3, one defines the variational solution to \eqref{e1.1}  through an SVI and one formulates the main existence result which is proved in Section 5, via the mentioned scaling method. In Section 6,
we prove  the positivity of solutions with nonnegative initial data and, in Section 7, we prove  the  finite time extinction  of solutions.

 We close this section with some remarks on our conditions (H1), (H2) and the stochastic integral in \eqref{e1.1}.

\begin{remark}\label{r11} {\rm \
\begin{itemize}\item[(i)] It is easy to check that under (H1) the sum in  \eqref{e1.2} converges in $L^2(\ooo;C([0,T];C(\ov\calo)))$  and that under (H1) and (H2) the sum in \eqref{e1.2} converges  in $L^2(\ooo;C([0,T];C^1(\ov\calo))).$ In particular, for $\mathbb{P}$-a.e. $\oo\in\ooo$ the map $$[0,T]\times\ov\calo\ni(t,\xi)\longmapsto W(t,\xi) (\oo)\in\rr$$is \ct\ and,  for each  $\xi\in\calo$, the process $(W(t,\xi))_{t\ge0}$ is  a real-valued (not standard) $(\calf_t)$-Brownian motion with quadratic variation $\mu(\xi)t$, with $\mu$ as defined in \eqref{12prim}. Furthermore, by Fernique's theorem,
\begin{equation}\label{e1.2prim}
\exp\left(\sup\limits_{0\le t\le T}|W(t)|_\9\right)\in L^p(\ooo)\mbox{\ \ for all }p\in(0,\9).\end{equation}

\item[(ii)] Let $F:[0,T]\times\calo\times\ooo\to\rr$ be such that $F$ (restricted to $[0,t])$ is $\mathcal{B}([0,t])\otimes\mathcal{B}(\calo)\otimes\calf_t$--measurable and $F\in\mathcal{L}^2([0,T]\times\calo\times\ooo)$, where $\mathcal{L}^2$ (instead of $L^2)$ denotes square integrable functions (rather than equivalence classes thereof). Then, for $\xi\in\calo$, we have $\pas$
    \begin{equation}\label{1.2}\int^t_0F(s,\xi)dW(s,\xi)
    =\dd\sum^\9_{k=1}\mu_k e_k(\xi)\int^t_0 F(s,\xi)d\b_k(s),\ t\in[0,T],\end{equation}
    where the sum on the right hand side converges in $L^2(\ooo;C([0,T];\rr))$ for each $\xi\in\calo$ and also in $L^2(\ooo;C([0,T];L^2(\calo)))$. Indeed, defining for $N\in\nn$
    $$W_N(t,\xi):=\dd\sum^N_{k=1}\mu_k e_k(\xi)\b_k(t),\ t\in[0,T],\ \xi\in\calo,$$
    we have for fixed $\xi\in\calo$, $N\in\nn$, by Doob's inequality and It\^o's isometry
    $$\barr{l}
    \E\left[\dd\sup_{t\in[0,T]}\left|\dd\int^t_0 F(s,\xi)d(W-W_N)(s,\xi)\right|^2\right]\vsp\qquad
    \le2\E\left[\dd\int^T_0|F(s,\xi)|^2ds\right]
    \dd\sum^\9_{k=N+1}\mu^2_k e^2_k(\xi)\earr$$and, similarly, for $N<M$,
    $$\barr{l}
    \E\left[\dd\sup_{t\in[0,T]}\int_\calo\left|
    \dd\sum^M_{k=N}\mu_ke_k(\xi)
    \dd\int^t_0 F(s,\xi)\b_k(ds)\right|^2d\xi\right]\vsp
    \qquad\le2\dd\sum^M_{k,k'=N}\mu_k\mu_{k'}\int_\calo e_k(\xi)e_{k'}(\xi)A_{k,k'}(\xi)d\xi,\earr$$where
    $$\barr{l}
    A_{k,k'}(\xi):=
    \E\left[\dd\int^T_0F(s,\xi)d\b_k(s)\dd\int^T_0 F(s,\xi)d\b_{k'}(s)\right]\vsp
    \qquad=\de_{k,k'}\E\left[\dd\int^T_0|F(s,\xi)|^2ds\right]\earr$$
    since  $\b_k,\b_{k'}$ are independent. So, both claimed convergences follow from (H1).

   \item[(iii)]  The assumption that $e_k$, $k\in\nn$, is an eigenbasis of the Dirichlet Laplacian is only used in the proof of Proposition \ref{l4.2} (ii) below. As it is pointed out there, this assumption is not necessary, provided the initial condition $x$ is in $H^1_0(\calo)$. Since Proposition \ref{l4.2} (ii) is only used in this paper for $x\in H^1_0(\calo)$, we may drop the above assumption on $e_k$, $k\in\nn$, and just assume that it is any orthonormal basis of $L^2(\calo)$ in $C^2(\ov\calo)\cap H^1_0(\calo)$.
\end{itemize}}\end{remark}

\begin{remark}\label{r1.2} {\rm Let $H=L^2(\calo)$ with usual inner product $\<\cdot,\cdot\>$ and norm $|\cdot|_2$. Under (H1), for all
$h\in H$  we have
$$\sum^\9_{k=1}|\mu_k|\,|\!\<h,e_k\>\!|\,|e_k|_\9\le C_\9|h|_2.$$Hence
$$\sum^\9_{k=1}\mu_k\<h,e_k\>e_k\in C(\ov\calo)$$and, for every $x\in H$, the following operator is well-defined
$$B(x)h:=x\dd\sum^\9_{k=1}\mu_k\<h,e_k\>e_k\ \left(=\dd\sum^\9_{k=1}\mu_k\<h,e_k\>(e_k\cdot x)\right),\ h\in H.$$
It is easy to check that $B(x)\in L_2(H,H)$ (= all Hilbert-Schmidt operators from $H$ to $H$) and that
\begin{equation}\label{e1.2primm}
\|B(x)\|_{L_2(H,H)}=\(\sum^\9_{k=1}\mu^2_k|e_kx|^2_2\)^{1/2}\le C_\9|x|_2.\end{equation}Therefore, if we consider the cylindrical Wiener process
$$\wt W(t):=(\b_k(t)e_k)_{k\in\mathbb{N}},$$
then, it is easy to check that, if $F$ is as in Remark \ref{r11} (ii), hence $(s,\oo)\mapsto F(s,\cdot,\oo)\in L^2(\calo)$ progressively measurable, then we have the following identities of $L^2(\calo)$-valued stochastic integrals in $L^2(\ooo;C([0,T];L^2(\calo)))$
 \begin{equation}\label{e1.3prim}
 \dd\int^\bullet_0 F(s)dW(s)=\int^\bullet_0 B(X(s))d\wt W(s)=\sum^\9_{k=1}\mu_k\int^\bullet_0F(s)e_k d\b_k(s),\end{equation}
 where also the sum on the right hand side converges in $L^2(\ooo;C([0,T];L^2(\calo)))$.  In particular, the stochastic integral in \eqref{e1.1} is a standard one.
  An easy application of the stochastic Fubini Theorem (cf. the proof of Claim 2 in the proof of Proposition \ref{p4.3}) then shows that, by \eqref{e1.3prim} and Remark \ref{r11} (ii), for every $t\in[0,T]$ and $\mathbb{P}$-a.e. $\oo\in\ooo$,
 $$\xi\mapsto \int^t_0 F(s,\xi) dW(s,\xi)(\oo)$$
 (which is a real-valued stochastic integral) is a $d\xi$-version of $\dd\int^t_0 F(s)dW(s)$ (which is an $L^2(\calo)$-valued stochastic integral).}\end{remark}


\section{Definition of stochastic variational solutions\\ and the main existence result}
\setcounter{equation}{0}

\begin{definition}\label{d3.1} {\rm Let $0<T<\9$ and let $x\in L^2(\calo)$. A stochastic process $X:[0,T]\times\Omega \to L^2(\calo)$ is said to be a variational solution to \eqref{e1.1} if the \fwg\ conditions hold.
\bit\item[(i)] $X$ is $(\calf_t)$-adapted,  has $\pas$  \ct\ sample paths in $L^2(\calo)$ and $X(0)=x$.
\item[(ii)] $X\in  L^2([0,T]\times\ooo;L^2(\calo)),$\, $\phi(X)\in L^1([0,T]\times\Omega).$
\item[(iii)] For each $(\calf_t)$- progressively measurable  process  $G\in L^2([0,T]\times\ooo;L^2(\calo))$ and each $(\calf_t)$-adapted  $L^2(\calo) $-valued process $Z$ with $\pas$ continuous sample paths such that $Z\in L^2([0,T]\times\ooo;H^1_0(\calo))$ and,  solving the equation\vspace*{-6mm}\eit
\begin{equation}\label{e3.1}
Z(t)-Z(0)+\int^t_0G(s)ds=\int^t_0Z(s)dW(s),\  \ t\in[0,T], \end{equation}
we have

\begin{equation}\label{e3.2}
\barr{l}
\dd\frac12\,\E|X(t)-Z(t)|^2_2+\E\dd\int^t_0\phi(X(\tau))d\tau
\le\dd\frac12\,\E|x-Z(0)|^2_2\vsp
\qquad+\E\dd\int^t_0\phi(Z(\tau))d\tau +\dd\frac12\,\E\int^t_0\int_\calo\mu(X(\tau)-Z(\tau))^2d\xi\,d\tau\vsp
\qquad+\E\dd\int^t_0
\<X(\tau)-Z(\tau),G(\tau)\>d\tau,\   t\in[0,T].\earr\end{equation}
Here, $\phi$ is defined by \eqref{e2.4},  $\mu=\dd\sum^\9_{k=1}\mu^2_ke^2_k$ and $\<\cdot,\cdot\>$ is the duality pairing   with pivot space $L^2(\calo)$. We also recall that \eqref{e3.1} has a unique solution for a given initial condition in $L^2(\calo)$.}\end{definition}

The relationship between \eqref{e1.1} and \eqref{e3.2} becomes more transparent if we recall that \eqref{e1.1} can be rewritten as \eqref{e2.9} and so we have
\begin{equation}\label{e3.3}
d(X-Z)+(\pp\phi(X)-G)dt\ni (X-Z)dW.\end{equation}
If we (formally) apply the It\^o  formula to $\dd\frac12\,|X-Z|^2_2$ in \eqref{e3.3} and take into account \eqref{e2.5}, we obtain just \eqref{e3.2} after taking expectation.
It should be emphasized, however, that $X$ arising in De\-fi\-ni\-tion \ref{d3.1} is not a strong solution to \eqref{e1.1} (or~\eqref{e2.9}) in the standard sense, that is,
$$X(t)-x\in-\int^t_0\pp\phi(X(s))ds+\int^t_0X(s)dW(s),\ \ff t\in(0,T).$$
We also note   that this concept of solution for \eq\ \eqref{e1.1} was already introduced in \cite{3}.
Theorem \ref{t3.1} below is our first main result.

\begin{theorem}\label{t3.1}  Let $\calo$ be a bounded and convex open subset of $\rr^N$ with smooth boundary.  For each $x\in  L^2(\calo) $ there is a  variational solution $X$ to \eq\ \eqref{e1.1}, such that, for all $p\in [2,\9)$,
\begin{equation}\label{e33prim}
\sup_{t\in[0,T]}\E[|X(t)|^p_2]\le\exp\left[C^2_\9\,\frac p2\,(p-1)\right]\|x\|^p_2.\end{equation}
$X$ is the unique solution in the class of all solutions $X$ such that, for some $\de>0$,

\begin{equation}\label{e34}
X\in L^{2+\de}(\ooo;L^2([0,T];L^2(\calo))).\end{equation}
Furthermore, if $x,x^*\in L^2(\calo)$ and $X,X^*$ are the corresponding variational solutions with initial conditions $x,x^*$, respectively, then, for some positive constant $C=C(N,C^2_\9)$,
\begin{equation}\label{e35}\E\left[\sup_{\tau\in[0,T]}|X(\tau)-X^*(\tau)|^2_2\right]\le 2|x-x^*|^2_2 e^{CT}.\end{equation}
In particular,
\begin{equation}\label{e35prim}
\E\left[\sup_{t\in[0,T]}|X(t)|^2_2\right]\le 2|x|^2_2 e^{CT},\end{equation}
and, moreover,    $X\in L^2(\ooo;C([0,T];L^2(\calo))).$
\end{theorem}

\begin{remark}\label{r3.1} {\rm A similar result was established in \cite{3} for the equation with additive noise where $N=1,2.$  However, in the definition of the solution in \cite{3}, erroneously was taken  the functional $\phi_0$  instead of $\phi$ defined above,  as it is correct. In this context we cite also the work \cite{6aa}, where this point was already clarified. Furthermore, by Remark \ref{r9.1} below, the convexity assumption on $\calo$ can be relaxed. It is enough that $\pp\calo$ can be parametrized locally by a convex $C^2$-map.}
\end{remark}

\begin{remark}\label{r3.2}{\rm  Apparently, the variational solution $X$ defined above does not satisfy in any common sense the Dirichlet homogeneous condition on $\partial\calo$, as written in \eqref{e1.1}. However, since $\phi(X)\in L^1((0,T)\times\Omega)$ and, as seen earlier,  $\phi$ is just  the closure in $L^2(\calo)$ of $\phi_0$ and, in particular, of the norm $|\na u|_1$ of the space $W^{1,1}_0(\calo),$  we may regard  $X$
as a generalized solution to \eqref{e1.1}. For instance, if in \eqref{e1.1} we replace the Dirichlet condition by the Neumann homogeneous condition, then in the above definition of the variational solution one should replace the function $\phi$~by
$$\phi_1(u) = ||Du||,\,  \forall  u\in BV(\calo)\cap L^2(\calo) ; \  \phi_1(u) = + \infty\  \mbox{otherwise}.$$}\end{remark}

\begin{remark}\label{r3.5}{\rm It follows from Lemma \ref{l7.3} below by Fatou's Lemma  that, for $N\le3$, in addition to \eqref{e35} we also have, for some constant $C>0$,
\begin{equation}\label{e3.7}
\E\left[\sup_{\tau\in[0,T]}|X(\tau)-X^*(\tau)|^N_N\right]\le 2|x-x^*|^N_Ne^{CT}.\end{equation}}\end{remark}

\begin{remark}\label{r3.6}{\rm
If $\in H^1_0(\calo)$, it follows by Lemma \ref{l5.1nou} and Fatou's Lemma that, for some $C>0$,
$$\E\left[\sup_{t\in[0,T]}\|X(t)\|^2_1\right]\le C\|x\|^2_1,$$hence $X\in L^2(\ooo;L^\9([0,T];H^1_0(\calo)))$. From this,  one can deduce that, if the initial condition $x$ is in $H^1_0(\calo)$, then the corresponding solution $X$ in Theorem \ref{t3.1} is, in fact, an ordinary variational solution of the (multivalued) equation \eqref{e1.1} (not just in the sense of SVI as in Definition \ref{d3.1}). Our main point is, however, here to have existence and uniqueness for all starting points $x\in L^2(\calo)$. Therefore, we skip the details on the simpler case of special initial conditions in $H^1_0(\calo)$.}\end{remark}

\section{The random differential equation\\ equivalent to \eqref{e1.1}}
\setcounter{equation}{0}

Inspired by \cite[Section 4]{3aa} and \cite{4}, we would like to   reduce \eq\ \eqref{e1.1} to the random dif\-fe\-ren\-tial \eq\

\begin{equation}\label{e4.1}
\barr{l}
\dd\frac{\pp Y}{\pp t}=e^{-W}{\rm div}({\rm sgn}(\na(e^WY)))-\dd\frac12\,\mu Y,\ \pas\ \mbox{in }(0,T)\times\calo,\vsp
Y(0,\xi)=x(\xi),\ \ \xi\in\calo,\vsp
Y=0\mbox{ on }(0,T)\times\pp\calo,
\earr\end{equation}by the substitution $Y(t)=e^{-W(t)}X(t)$. The meaning of  boundary condition in \eqref{e4.1} is taken in the generalized sense as   discussed in Remark \ref{r3.2}.
(We~note that, in equation \eqref{e1.1}, $XdW$ is meant to be an It\^o differential,    other\-wise, i.e., if it is taken in the Stratonovich sense,  then, in the corresponding equation \eqref{e4.1}, the linear term $\dd\frac12\,\mu Y$  would be
 missing.)

To do the reduction from \eqref{e1.1} to \eqref{e4.1} rigorously, our definitions of solutions for \eqref{e4.1} must be again in the sense of a variational inequality, but this time a deterministic one, since the test processes $\wt Z$ (replacing $Z$ in Definition \ref{d3.1}) solve a deterministic PDE, however, with random coefficients.

\begin{definition}\label{d4.1} {\rm Let $0<T<\9$ and let $x\in L^2(\calo)$. A stochastic process  $Y:[0,T]\times\ooo\to L^2(\calo) $ is said to be a {\it variational solution} to \eqref{e4.1} if the \fwg\ conditions hold:
\bit
\item[(i)] $Y$ is    $(\calf_t)$-adapted, has $\pas$ \ct\ sample paths, and $Y(0)=x.$
\item[(ii)] $e^WY \in L^2([0,T]\times\ooo;L^2(\calo)), \, \phi(e^{W}Y)\in L^1([0,T]\times\Omega).$
\item[(iii)] For each   $(\calf_t)$-progressively measurable process $G\in L^2([0,T]\times\ooo;L^2(\calo))$ and each  $(\calf_t)$-adapted, $L^2(\calo)$-valued process $\wt Z$ with $\pas$ \ct\ sample paths such that $e^W\wt Z\in L^2([0,T]\times\ooo;H^1_0(\calo))$ and  solving    the \eq\vspace*{-4mm}\eit
\begin{equation}\label{e4.3a}
\barr{r}
\dd\wt Z(t)-\wt Z(0)+\int^t_0 e^{-W(s)}G(s)ds+\frac12\int^t_0\mu\wt Z(s)ds=0,\vsp     t\in[0,T],\ \pas,\earr\end{equation}
we have

\begin{equation}\label{e4.4}
\barr{l}
\dd\frac12\,\mathbb{E}|e^{W(t)}(Y(t)-\wt Z(t))|^2_2+\mathbb{E}\dd\int^t_0\phi(e^{W(\tau)}Y(\tau))d\tau\vsp
\qquad\le\dd\frac12\,\mathbb{E}|x-\wt Z(0)|^2_2 +\mathbb{E}\dd\int^t_0\phi(e^{W(\tau)}\wt Z(\tau))d\tau\vsp
\qquad\dd+\dd\frac12\,\mathbb{E}\dd\int^t_0\int_\calo\mu e^{2W(\tau)}(Y(\tau)-\wt Z(\tau))^2d\xi\,d\tau\vsp
\qquad+\mathbb{E}\dd\int^t_0{} \<e^{W(\tau)}(Y(\tau)-\wt Z(\tau)),G(\tau)\> d\tau,\  t\in[0,T].\earr\end{equation}}
\end{definition}
We recall that the deterministic equation \eqref{e4.3a} has a unique solution for a given initial condition in $L^2(\calo)$ for $\mathbb{P}$-a.e. given $\oo\in\ooo$.


\begin{proposition}\label{p4.2} $X:[0,T]\times\ooo\to L^2(\calo)$ is a variational solution to equation  \eqref{e1.1} if and only if $Y:=e^{-W}X$ is a variational solution to \eqref{e4.1}.
\end{proposition}

The above proposition is an immediate consequence of Proposition \ref{p4.3}(iii)  below, which addresses a technical, but very important issue. To be precise (and make our point) in its proof, we have to distinguish between the space $\mathcal{L}^2(\calo)$ of square integrable functions and $L^2(\calo)$, i.e., the corresponding $d\xi$-classes.

\begin{proposition}\label{p4.3} Let $G\in L^2([0,T]\times\ooo;L^2(\calo))$ be $(\calf_t)$-progressively measurable and $Z(0)\in L^2(\ooo,\calf_0;L^2(\calo)).$ Let $G^0$ be a $(dt\otimes d\xi\otimes \mathbb{P})$-version of $G$ such that $(t,\oo)\longmapsto G^0(t,\xi,\oo)$ is $(\calf_t)$-progressively measurable and in $L^2([0,T]\times\ooo)$ for every $\xi\in\calo$. Furthermore, let $Z^0(0)$ be a $(d\xi\otimes\mathbb{P})$-version of $Z(0)$ such that $\oo\longmapsto Z^0(0)(\xi,\oo)$ is $\calf_0$-measurable for all $\xi\in\calo$.
\bit\item[\rm(i)] Fix $\xi\in\calo$. Then
\begin{equation}\label{4.3a}
\barr{l}
Z^0_\xi(t):=e^{W(t,\xi)-\frac12\,\mu(\xi)t}Z^0(0)(\xi)\\
\qquad-e^{W(t,\xi)-\frac12\,\mu(\xi)t}\dd\int^t_0 e^{-W(s,\xi)+\frac12\,\mu(\xi)s}
G^0(s,\xi)ds,\  t\in[0,T],\earr\end{equation}
is a real-valued \ct\ solution to the stochastic differential equation
\begin{equation}\label{e4.6a}
\barr{l}
dZ^0_\xi(t)=-G^0(t,\xi)dt+Z^0_\xi(t)dW(t,\xi),\ t\in[0,T],\vsp Z^0_\xi(0)=Z^0(0,\xi),\earr\end{equation}
which is $\mathcal{B}([0,t])\otimes\mathcal{B}(\calo)\otimes\calf_t$-measurable for each $t\in[0,T].$

Furthermore, the map $[0,T]\ni t\longmapsto Z^0_\cdot(t)\in\mathcal{L}^2(\calo)$ is $\pas$ con\-ti\-nous. Hence the corresponding $d\xi$-classes $Z(t)\in L^2(\calo)$, $t\in[0,T]$, form the unique solution to \eqref{e3.1}.

\item[\rm(ii)] Fix $\xi\in\calo$. Then
$$\barr{r}
\wt Z^0_\xi(t):= e^{-\frac12\,\mu(\xi)t}Z^0(0)(\xi)-e^{-\frac12\,\mu(\xi)t}\dd\int^t_0 e^{-W(s,\xi)+\frac12\,\mu(\xi)s}G^0(s,\xi)ds,\vsp t\in[0,T],\earr$$is a real-valued \ct\ solution to the differential equation
$$d\wt Z^0_\xi(t)=-e^{-W(t,\xi)}G^0(t,\xi)dt-\frac12\,\mu(\xi)\wt Z^0_\xi(t)dt,\ \wt Z^0_\xi(0)=Z^0(0,\xi),$$
which is $\mathcal{B}([0,t])\otimes\mathcal{B}(\calo)\otimes\calf_t$-measurable for each $t\in[0,T].$

Furthermore, the map
 $[0,T]\ni t\longmapsto\wt Z^0_\cdot(t)\in\mathcal{L}^2(\calo)$ is $\pas$ con\-ti\-nuous. Hence the corresponding $d\xi$-classes $\wt Z_\cdot(t)\in L^2(\calo)$, $t\in[0,T]$, form the unique solution of the deterministic equation \eqref{4.3a} for $\mathbb{P}$-a.e. given $\oo\in\ooo$.

\item[\rm(iii)] An $(\calf_t)$-adapted $\pas$ \ct\ $L^2(\calo)$-valued process $(Z(t))_{t\in[0,T]}$ is a solution to the stochastic equation \eqref{e3.1} if and only if $(e^{-W(t)}Z(t))_{t\in[0,T]}$ is a solution to the deterministic equation \eqref{e4.1} for $\mathbb{P}$-a.e. given $\oo\in\ooo$.
\eit
\end{proposition}

\n{\bf Proof.} (iii) is an immediate consequence of (i) and (ii). (ii) is more or less well-known since it is about a deterministic equation and the proof is anyway similar to that of (i). Therefore, we only prove (i).

First, we note that applying a mollifier in $\xi$ and taking the limsup of a properly chosen subsequence, the mentioned version of $G^0$ and $Z^0(0)$ always exist.
Obviously, $Z^0_\xi$ is a well-defined, $(\calf_t)$-adapted, $\pas$ \ct\ real-valued process, and applying It\^o's product formula we obtain that it solves \eqref{e4.6a}. Furthermore, the stated continuity in $\mathcal{L}^2(\calo)$ is obvious. So, it remains to show the last part of the assertion, which follows from the following two claims.

\mk\n{\bf Claim 1.} Let $t\in[0,T].$ Then $\mathbb{P}$-a.e.
 $\xi\longmapsto \int^t_0 G^0(s,\xi)ds,\ \ \xi\in\calo,$ is a $d\xi$-version of the $L^2(\calo)$-valued Bochner integral $\int^t_0G(s)ds.$

\mk\n{\bf Claim 2.} Let $t\in[0,T].$ Then $\mathbb{P}$-a.s.
 $\xi\longmapsto \int^t_0 Z^0_\xi(s)dW(s,\xi),\ \ \xi\in\calo,$
is a $d\xi$-version of the $L^2(\calo)$-valued stochastic  integral $\int^t_0Z(s)dW(s).$\mk

Claim 1 is a trivial consequence of Fubini's theorem. So, we only prove Claim 2 whose proof is similar, but based on the stochastic Fubini theorem.

\mk\n{\bf Proof of Claim 2.} Let $i\in\mathbb{N}.$ Then $\mathbb{P}$-a.s. for every $t\in[0,T]$, by Remark \ref{r11} (ii),
$$\barr{l}
\dd\int_\calo e_i(\xi)\dd\int^t_0 Z^0_\xi(s)dW(s,\xi)d\xi\vsp
\qquad=  \dd\sum^\9_{k=1}\mu_k\dd\int_\calo e_i(\xi)e_k(\xi)\dd\int^t_0Z^0_\xi(s)d\beta_k(s)d\xi\vsp
\qquad=\dd  \sum^\9_{k=1}\mu_k\dd\int^t_0
  \dd\int_\calo
e_i(\xi)e_k(\xi) Z^0_\xi(s) d\xi d\beta_k(s)\vsp
\qquad= \dd\sum^\9_{k=1}\mu_k\dd\int^t_0
\<e_i,e_k Z(s)\>d\beta_k(s)\vsp
\qquad=  \dd\sum^\9_{k=1}\mu_k\<e_i,\dd\int^t_0 e_k Z(s)
d\beta_k(s)\>\vsp
\qquad=\<e_i,\dd\int^t_0 Z(s)dW(s)\>,\earr$$where we used the stochastic Fubini theorem in the second equality. Now, Claim 2 follows. $\Box$

\begin{remark}\label{r4.4} {\rm Proposition \ref{p4.3} justifies to apply It\^o's formula for a solution $Z(t),$ $t\in[0,T]$, to \eqref{e3.1} for $d\xi$-a.e. $\xi\in\calo$ to the process $Z(t)(\xi)$, $t\in[0,T]$, by taking the version $Z^0_\xi(t)$, $t\in[0,T]$, from Proposition \ref{p4.3}(i).    We stress that for Proposition \ref{p4.3} we only used (H1), not (H2) (see Remark \ref{r11}).}\end{remark}

In particular, by Theorem \ref{t3.1},   Proposition \ref{p4.2} and \eqref{e59prim} below, we have the following existence result for \eqref{e4.1}, which has an intrinsic interest.

\begin{proposition}\label{p4.5}  Under the assumptions of Theorem {\rm\ref{t3.1}},  for each $x\in L^2(\calo)$, there is a  variational solution $Y$ to \eqref{e4.1}, which is unique in the class of all solutions $Y$ such that, for some $\de>0$,
$$Y\in L^{2+\de}(\ooo;L^2([0,T];L^2(\calo))).$$
Moreover, $$Y\in L^2(\ooo; C([0,T];L^2(\calo))).$$
\end{proposition}

\section{Proof of Theorem \ref{t3.1}}
\setcounter{equation}{0}

It should be said that, for the proof of the uniqueness part of Theorem \ref{t3.1}, as well as for the finite-time extinction property of the solutions to \eqref{e1.1}, it is convenient and apparently necessary to replace \eqref{e1.1} by \eqref{e4.1} and to construct approximating schemes for both equations.

We approximate \eqref{e1.1} by
\begin{equation}\label{e5.1a}\barr{l}
dX_\lbb={\rm div}\ \wt\psi_\lbb(\na X_\lbb)dt+X_\lbb dW\ \inr\ (0,T)\times\calo,\vsp
X_\lbb(0)=x\ \inr\ \calo,\vsp
X_\lbb=0\ \onr\ (0,T)\times\pp\calo,\earr\end{equation}
and the corresponding rescaled  equation \eqref{e4.1} by
\begin{equation}\label{e4.5}
\barr{l}
\dd\frac{dY_\lbb}{dt}=e^{-W}{\rm div}(\wt\psi_\lbb(\na(e^WY_\lbb)))-\dd\frac12\,\mu Y_\lbb\vsp\hfill \mbox{in }(0,T)\times\calo,\ \pas,\vsp
Y_\lbb(0)=x\ \mbox{in }\calo,\ \ Y_\lbb=0\mbox{ on }(0,T)\times\pp\calo,\earr\end{equation}
where $\lbb\in(0,1]$,  $\wt\psi_\lbb(u) = \psi_\lbb(u) +\lbb u ,\ \forall u\in  \rr^{N}$.

 In \eqref{e4.5}, $\frac d{dt}\,Y_\lbb\in L^2(0,T;H\1(\calo))$ is the strong derivative of $t\to Y_\lbb(t)$ and the operator div is taken in sense of distributions on $\calo$.

  Here, $\psi_\lbb$ is the Yosida approximation of the   function $\psi(u)={\rm sgn}\ u$, that is (see, e.g.,~\cite{2}),
\begin{equation}\label{e4.6}
\psi_\lbb(u)=\left\{\barr{lll}
\dd\frac1\lbb\,u&&\mbox{if }|u|\le\lbb,\vsp
\dd\frac u{|u|}&&\mbox{if }|u|>\lbb.\earr\right.\end{equation}


Let  $j_\lbb(u)=\dd\inf_v\left\{\dd\frac{|u-v|^2}{2\lbb}
+|v|\right\}$  be the Moreau--Yosida approximation of the function $v\to |v|.$  We  recall that $\na j_\lbb=\psi_\lbb$, $\ff \lbb>0$ (see, e.g., \cite{2}, p.~48).  We first prove   the existence of a strong
solution $Y_\lbb$ to \eqref{e4.5}.

It should be emphasized that, for the existence and uniqueness part of the proof, it is convenient to analyze equation \eqref{e5.1a} while, for getting sharp estimates on the variational solutions $X$ to \eqref{e1.1}, it is necessary to work directly with the random equation \eqref{e4.5} instead of \eqref{e5.1a}. As regards the existence and uniqueness for \eqref{e5.1a}, \eqref{e4.5}, we have:

\begin{proposition}\label{l4.2}\
\begin{itemize}
\item[\rm(i)] For each $\lbb\in(0,1]$ and each $x\in L^2(\calo)$, there is a unique strong solution $X_\lbb$ to \eqref{e5.1a} which sa\-tisfies $X_\lbb(0)=x$, that is, $X_\lbb$ is $\pas$ con\-ti\-nuous in $L^2(\calo)$ and $\{\calf_t\}$-adapted such that
\begin{equation}\label{e5.4a}
\barr{l}
\dd X_\lbb\in L^2([0,T]\times\ooo;H^1_0(\calo)),\vsp
X_\lbb(t)=x+\dd\int^t_0{\rm div}\,\wt\psi_\lbb
\(\na X_\lbb(s)\)ds+\dd\int^t_0 X_\lbb(s)dW(s),\vsp\hfill t\in[0,T],\ \pas\earr\end{equation}
Furthermore, $X_\lbb\in L^2(\ooo;C([0,T];L^2(\calo)))$ and, for all $p\in [2,\9),$
\begin{equation}\label{e54prima}
\sup_{t\in[0,T]}\E\left[|X_\lbb(t)|^p_2\right]\le
\exp\left[C^2_\9\,\frac p2\,(p-1)\right]|x|^p_2,\end{equation}
and,
if $x,x^*\in  L^2(\calo)$ and $X_\lbb$ and $X^*_\lbb$ are the corresponding solutions with initial conditions $x,x^*$, respectively, then, for some positive constant $C=C(N,C^2_\9)$,
\begin{equation}\label{e54prim}
\E\left[\sup_{\tau\in[0,T]}|X_\lbb(\tau)-X^*_\lbb(\tau)|^2_2\right]
\le2|x-x^*|^2_2e^{CT}.\end{equation}

\item[\rm(ii)] $Y_\lbb=e^{-W}X_\lbb$ is an  $(\calf_t)$-adapted process $Y_\lbb:[0,T]\times\ooo\to L^2(\calo)$  with $\pas$ continuous paths which is the unique solution of \eqref{e4.5}, i.e., it satisfies $\pas$ equation \eqref{e4.5} with $Y_\lbb(0)=x$ and
\begin{equation}\label{e4.7}
 Y_\lbb\in L^2([0,T];H^1_0(\calo))\cap C([0,T]; L^2(\calo))\cap W^{1,2}([0,T];H\1(\calo)),
\end{equation}
 a.e. $t\in[0,T]$.
 \item[\rm(iii)]  If $x\in H^1_0(\calo)$, then $\pas$
\begin{equation}\label{5.5a}
X_\lbb\in C([0,T];H^1_0(\calo))\end{equation}and
\begin{equation}\label{e57prim}
X_\lbb\in L^2([0,T]\times\ooo;H^2(\calo)).\end{equation}\end{itemize}
\end{proposition}

\begin{remark}\label{r51prim} {\rm
It is readily seen that, by It\^o's formula, $X_\lbb$ is also a variational solution to \eqref{e5.1a} in the sense of Definition \ref{d3.1}, where $\phi(y)$ is replaced by
$$ \wt\phi_\lbb(y)= \int_\calo\(j_\lbb(\na y)+\frac\lbb2\,|\na y|^2\)d\xi.$$
}\end{remark}

 \n{\bf Proof of Proposition \ref{l4.2}.} Consider the operator $A_\lbb:H^1_0(\calo)\to H\1(\calo)$ defined by
\begin{equation}\label{e5.6prim}
\<A_\lbb y,\vf\>=\int_\calo\wt\psi_\lbb(\na y)\cdot\na\vf d\xi,\ \ff \vf\in H^1_0(\calo),\end{equation}
and note that $A_\lbb$ is demicontinuous   (see, for instance, \cite{2}, p.~81).

Moreover, we have
$$\barr{rcll}
\|A_\lbb y\|_{-1}&\le&\lbb\|y\|_1+\(\dd\int_\calo d\xi\)^{\frac12},&\ff y\in H^1_0(\calo),\vsp
\<A_\lbb y_1-A_\lbb y_2,y_1-y_2\>&\ge&\lbb\|y_1-y_2\|^2_1,&\ff y\in H^1_0(\calo).\earr$$
On the other hand, equation \eqref{e5.1a} can be rewritten as
\begin{equation}\label{e5.7a}
\barr{l}
dX_\lbb+A_\lbb X_\lbb dt=X_\lbb dW,\ t\in[0,T],\vsp
X_\lbb(0)=x.\earr\end{equation}
Then, by the standard existence theory for stochastic differential equations associated with nonlinear monotone and demicontinuous operators in a dua\-lity pair $(V,V')$ (\cite{18a}, \cite{21a}, \cite{11}) equation \eqref{e5.7a} (equivalently, \eqref{e5.1a}) has a unique strong solution
$X_\lbb$ satisfying \eqref{e5.4a} and \eqref{e54prim}. \eqref{e54prima} is then an easy consequence of It\^o's formula for $|X_\lbb|^2_2$ (see, e.g., \cite{11}).

 To  prove (ii), below we use $\<\cdot,\cdot\>_2$ to denote the inner product in $L^2(\calo)$, in order to avoid confusion with the quadratic variation process.

Let $\vf\in H^1_0(\calo)\cap L^\9(\calo)$. Then, for every $t\in[0,T]$,
$$\<\vf,e^{-W(t)}X_\lbb(t)\>_2
=\dd\sum^\9_{j=1}\<e_j,e^{-W(t)}\vf\>_2\<e_j,X_\lbb(t)\>_2.$$
Furthermore, by It\^o's formula and Remark \ref{r11}, we have $d\xi\otimes\mathbb{P}$-a.e.  that, for all $t\in[0,T]$,
$$e^{-W(t,\xi)}=1-\dd\int^t_0e^{-W(s,\xi)}dW(s,\xi)
+\dd\frac12\,\mu(\xi)\dd\int^t_0e^{-W(s,\xi)}ds.$$
Now, fix $j\in\nn$. Then, by Remark \ref{r11} (ii) and Remark \ref{r1.2}, we have $\mathbb{P}$-a.e. that, for all $t\in[0,T]$,
$$\barr{lcl}
\<e_j,e^{-W(t)}\vf\>_2&=&
\<e_j,\vf\>_2-\dd\sum^\9_{k=1}\mu_k
\dd\int_\calo e_j(\xi)\vf(\xi)e_k(\xi)
\dd\int^t_0 e^{-W(s,\xi)}d\b_k(s)d\xi\vsp
&&+\,\dd\frac12\int^t_0\<e_j,\mu e^{-W(s)}\vf\>_2ds\vsp
&=&\<e_j,\vf\>-\dd\sum^\9_{k=1}\mu_k\dd\int^t_0
\<e_j,e_ke^{-W(s)}\vf\>_2d\b_k(s)
\vsp &&+\,\dd\frac12\int^t_0\<e_j,\mu e^{-W(s)}\vf\>_2ds,\earr$$
where we used the stochastic Fubini Theorem in the second equality and the sums converge in $L^2(\ooo;C([0,T];\rr)).$ By It\^o's product rule we hence obtain $\pas$ that, for all $t\in[0,T]$,
$$\barr{l}
\<e_j,e^{-W(t)}\vf\>_2\<e_j,X_\lbb(t)\>_2
=\<e_j,\vf\>_2\<e_j,x\>_2\vsp
\qquad+\dd\int^t_0\<e_j,e^{-W(s)}\vf\>_2\<e_j,{\rm div}\,\wt \psi_\lbb(\na X_\lbb(s))\>ds\vsp
\qquad+\dd\sum^\9_{k=1}\mu_k\dd\int^t_0
\<e_j,e^{-W(s)}\vf\>_2\<e_j,X_\lbb(s)e_k\>_2d\b_k(s)\vsp
\qquad\dd-\sum^\9_{k=1}\mu_k\dd\int^t_0
\<e_j,e_k e^{-W(s)}\vf\>_2\<e_j,X_\lbb(s)\>_2d\b_k(s)\earr$$
$$\barr{l}
\qquad+\dd\frac12\int^t_0\<e_j,\mu e^{-W(s)}\vf\>_2\<e_j,X_\lbb(s)\>_2ds\vsp
\qquad-\dd\sum^\9_{k=1}\mu^2_k\int^t_0
\<e_j,X_\lbb(s)e_k\>_2\<e_j,e_k
e^{-W(s)}\vf\>_2ds,\earr$$
where all the sums converge in $L^2(\ooo;C([0,T];\rr))$ and  interchanging the infinite sums with stochastic differentials is justified by Remark \ref{r11} (ii) and Remark \ref{r1.2}, because of \eqref{e54prima} and since, by \eqref{e1.2prim},
\begin{equation}\label{e59prim}
\sup_{(t,\xi)\in[0,T]\times\calo} e^{-W(s,\xi)}\
|X_\lbb|_2\in L^p([0,T]\times\ooo;\rr)\mbox{ for all }p\ge1.\end{equation}
(We shall implicitly use both \eqref{e54prima} and \eqref{e59prim}  several times in the rest of this paper without further notice.)

Now, we sum the above equation from $j=1$ to $j=\9$ and assume that we can interchange this summation both with the sum over $k$ and with  the deterministic and stochastic integrals (which we shall justify below). Then, because the two terms involving the stochastic integrals cancel, we obtain
$$\barr{l}
\<\vf,e^{-W(t)}X_\lbb(t)\>_2
=\<\vf,x\>_2+\dd\int^t_0\<\vf,e^{-W(s)}{\rm div}\,\wt\psi_\lbb
(\na X_\lbb(s))\>ds\vsp
\qquad+\dd\frac12\int^t_0\<\vf,\mu e^{-W(s)}X_\lbb(s)\>_2ds -\dd\sum^\9_{k=1}\mu^2_k\int^t_0\<\vf,e^2_k e^{-W(s)}X_\lbb(s)\>_2ds,\earr$$
which immediately implies that $Y_\lbb=e^{-W}X_\lbb$ solves \eqref{e4.5}.

To justify interchanging sums and integrals, it suffices to note that, for the second term on the right hand side, this is true because $\{e_k\}$ is the eigenbasis of the Laplacian and that for the last term this is obvious because of (H1), while, for the two terms  which cancel each other and involve stochastic integrals, this follows by applying the Burkholder--Davis--Gundy inequality and (H1). If, however, $x\in H^1_0(\calo)$, then, by Lemma \ref{l5.1nou} below, ${\rm div}\ \wt\psi_\lbb(\na X_\lbb)\in L^2([0,T]\times\ooo;L^2(\calo))$ (and not only in $L^2([0,T]\times\ooo;H^{-1}(\calo)))$. Hence, the above equality is true for any orthonormal basis $e_k$, $k\in\nn$, of $L^2(\calo)$  in $C^2(\ov\calo)$.


It remains to prove the uniqueness. In fact, as it will be explained below, by standard methods one can prove directly the existence and uniqueness of a solution $Y_\lbb$ to \eqref{e4.5}, which hence must be of the form $Y_\lbb= e^{-W}X_\lbb$. To this end, for each $\oo\in\ooo$, consider the operator $$\wt A=\wt A_\lbb(t,\oo):H^1_0(\calo)\to H\1(\calo)$$ defined by
\begin{equation}\label{5.7a}
\barr{l}
\<\wt A_\lbb(t)y,\vf\>=\dd\int_\calo\psi_\lbb(\na e^{W(t)}y))\cdot\na(e^{-W(t)}\vf)d\xi\vsp
\qquad+\lbb\dd\int_\calo\na(e^{W(t)}y)\cdot\na(e^{-W(t)}\vf)d\xi+
\dd\frac12\int_\calo\mu y\vf\,d\xi,\vsp\hfill \ff\vf\in H^1_0(\calo).\earr\end{equation}

In terms of $\wt A_\lbb$, equation \eqref{e4.5} becomes
\begin{equation}\label{e5.7b}
\barr{l}
\dd\frac{dY_\lbb}{dt}+\wt A_\lbb(t)Y_\lbb(t)=0,\mbox{ a.e. }t\in(0,T),\vsp
X_\lbb(0)=x.\earr\end{equation}
It is easily seen that, for every $t\in[0,T]$ and $\pas$, $\oo\in\ooo,$ $\wt A_\lbb(t)=\wt A_\lbb(t)(\oo)$ is demicontinuous (that is, strongly-weakly continuous), coercive, that is,
\begin{equation}\label{e5.8a}
\<\wt A_\lbb(t)y,y\>\ge\lbb\|y\|^2_1-\a^\lbb_t|y|^2_2,\ \ff y\in H^1_0(\calo),\end{equation}bounded, that is,
\begin{equation}\label{e511prim}
\|\wt A_\lbb(t)y\|_{-1}\le C_t(1+\|y\|_1),\ \ff y\in H^1_0(\calo),\end{equation}
and $\de$-monotone, that is,
\begin{equation}\label{e511primm}
\<\wt A_\lbb(t)y-\wt A_\lbb(t)z,y-z\>+\de^\lbb_t|y-z|^2_2\ge0,\ \ff y,z\in H^1_0(\calo),\end{equation}
where $C_t,\a^\lbb_t,\de^\lbb_t:\ooo\to\rr_+$, $t\in[0,T]$, are $(\calf_t)$-adapted processes, $\pas$ continuous   on $[0,T]$. (Since, as pointed out before, we only need the uniqueness part, i.e., we only need \eqref{e511primm}, for the reader's convenience we include its proof in Appendix 2, i.e., Section 10.)


Hence, for each $x\in L^2(\calo)$, there is a unique solution $Y_\lbb$ to \eqref{e5.7b} sa\-tis\-fying \eqref{e4.7}. (See, e.g., \cite{5}, p. 177). This completes the proof of (ii). To prove (iii), we need the following two lemmas:

\begin{lemma}\label{l5.1nou} Let $x\in H^1_0(\calo)$. Then,
$X_\lbb\in L^2(\ooo;L^\9([0,T];H^1_0(\calo)))\cap\break L^2([0,T]\times\ooo;H^2(\calo))$ and
\begin{equation}\label{5.12}
\E\left[\sup_{t\in[0,T]}\|X_\lbb(t)\|^2_1\right]+
\lbb\E\int^T_0|\D X_\lbb(t)|^2_2dt\le C\|x\|^2_1,\ \lbb\in(0,1].\end{equation}
\end{lemma}


 \n{\bf Proof of Lemma \ref{l5.1nou}.} In this proof, constants may change from line to line, though we continue to denote them by $C$. We set $A=-\D,$ $D(A)=H^1_0(\calo)\cap H^2(\calo)$, $J_\vp=(1+\vp A)\1,$ $A_\vp=AJ_\vp=\frac1\vp\,(I-J_\vp)$ and note that, by virtue of Corollary \ref{c8.8} in Appendix 1, we have
 \begin{equation}\label{ez5.13}
 \barr{ll}
 -\<A_\vp X_\lbb,{\rm div}\,\psi_\lbb(\na X_\lbb)\>
 \!\!\!&=\dd\frac1\vp\int_\calo(\na y-\na J_\vp(y))\cdot\psi_\lbb(\na y)d\xi\vsp
 &\ge\dd\frac1\vp\int_\calo(j_\lbb(\na y)-j_\lbb(\na J_\vp(y)))d\xi\ge0.\earr\end{equation}

 Now, we apply It\^o's formula to the function $\vf(x)=\frac12\,|A^{\frac12}_\vp x|^2_2.$ We have $D \vf=A_\vp$, and so we get by Hypotheses (H1) and (H2)  that
 \begin{equation}\label{ez5.14}
 \barr{l}
 \dd\frac12\, |A^{\frac12}_\vp X_\lbb(t)|^2_2
 +\lbb \dd\int^t_0\<A_\vp X_\lbb(s),A X_\lbb(s)\>ds\vsp
 \qquad- \dd\int^t_0\<A_\vp X_\lbb(s),\mbox{div}\,\psi_\lbb(\na X_\lbb(s))\>ds\vsp
 \qquad\le\dd\frac12\,|x|^2_2
 +C \dd\int^t_0\|X_\lbb(s)\|^2_1ds\vsp\qquad
 +\dd\int^t_0\<A_\vp X_\lbb(s),X_\lbb(s)dW(s)\>,\ t\in[0,T],\earr\end{equation}
 since $|A^{\frac12}_\vp x|_2\le\|x\|_1$, $\ff x\in H^1_0(\calo),\ \vp\in(0,1].$

 Now, keeping in mind that, for all $\vp>0$,
 $$\<A_\vp y,Ay\>\ge|A_\vp y|^2_2,\ \ff y\in H^1_0(\calo),$$
 and, taking into account \eqref{ez5.13}, we obtain by \eqref{ez5.14} that, for some $C>0$ independent of $\lbb$ and $\vp$,
 \begin{equation}\label{e520prim}
\barr{l}
 \dd |A^{\frac12}_\vp X_\lbb(t)|^2_2+\lbb
 \int^t_0|A_\vp X_\lbb(s)|^2_2 ds\le
 \dd\frac12\,|x|^2_2+C \dd\int^t_0\|X_\lbb(s)\|^2_1ds\vsp
 \hfill
 +\dd\int^t_0\<A_\vp X_\lbb(s),X_\lbb(s)dW(s)\>,\  t\in[0,T],\ \ff\lbb,\vp>0.\earr\end{equation}

 In particular, for all $t\in[0,T]$
 $$\barr{r}
 \dd\sup_{r\in[0,t]}|A^{\frac12}_\vp X_\lbb(r)|^2_2\le\|x\|^2_1+C
 \dd\int^t_0\sup_{r\in[0,s]}\|X_\lbb(r)\|^2_1ds\vsp
 +\dd\sup_{r\in[0,t]}\left|\dd\int^r_0\<A_\vp X_\lbb(s),X_\lbb(s)dW(s)\>\right|.\earr$$
 Hence, by the Burkholder-Davis-Gundy (for $p=1)$ and Gronwall's inequalities, we obtain that, for some $C>0$, independent of $\lbb$ and $\vp$,
 $$\E\left[\sup_{s\in[0,T]}|A^{\frac12}_\vp X_\lbb(s)|^2_2\right]\le
 2\|x\|^2_1e^{CT},\ \ff\lbb,\vp\in(0,1].$$
  Letting $\vp\to0$, we obtain
\begin{equation}\label{ez5.20prim}
\E\left[\sup_{t\in[0,T]}\|X_\lbb(t)\|^2_1\right]\le C\|x\|^2_1,\ \ff\lbb\in(0,1].\end{equation}
Hence, taking expectation in \eqref{e520prim} and letting $\vp\to0$,  we obtain
 $$\lbb\E\int^T_0|\D X_\lbb(s)|^2_2ds\le C\|x\|^2_1,\ \ff \lbb\in(0,1].$$
 This completes the proof of Lemma \ref{l5.1nou}. $\Box$


\begin{lemma}\label{l5.4nou} Let   $x\in H^1_0(\calo)$. Then,   $Y_\lbb\in C([0,T];H^1_0(\calo))\cap L^2([0,T];H^2(\calo))$, $\pas$
\end{lemma}

\n{\bf Proof.} We rewrite \eqref{e4.5} as the linear parabolic random equation
\begin{equation}\label{5.27}
\barr{l}
\dd\frac{\pp Y_\lbb}{\pp t}=\lbb\D Y_\lbb+f(t,\xi)\ \mbox{in }(0,T)\times\calo,\vsp
 Y_\lbb=0\mbox{ on }(0,T)\times\pp\calo,\vsp
 Y_\lbb(0,\xi)=x(\xi)\mbox{ in }\calo,\earr\end{equation}
 where
 $$\barr{r}
f(t,\xi)=e^{-W(t,\xi)}{\rm div}\,\psi_\lbb
(\na(e^{W(t,\xi)}Y_\lbb(t,\xi)))
-\dd\frac12\,\mu(\xi)Y_\lbb(t,\xi)\vsp
\qquad+2\lbb\na W(t,\xi)\cdot\na Y_\lbb(t,\xi)+\D W(t,\xi)Y_\lbb(t,\xi)+Y_\lbb|\na W|^2.\earr$$
Since, for $y\in H^1_0(\calo)\cap H^2(\calo)$,
\begin{equation}\label{e521prim}
\barr{l}
{\rm div}\,\psi_\lbb(\na y)=
\left\{\barr{ll}
\dd\frac1\lbb\,\D y&\mbox{ on } \{|\na y|\le\lbb\},\vsp
\dd\frac{\D y}{|\na y|}-\frac{\na y\cdot\na|\na y|}{|\na y|^2}&\mbox{ on }\{|\na y|>\lbb\},\earr\right.\earr\end{equation}
by Lemma \ref{l5.1nou}, we know that
\begin{equation}\label{5.28}
f(t)\in L^2(0,T;L^2(\calo)),\ \pas\end{equation}
Then, by the general theory of linear parabolic equations (see, e.g. \cite{H}), we have, for each $\oo\in\ooo$,
\begin{equation}\label{5.29}
Y_\lbb\in C([0,T];H^1_0(\calo))\cap L^2([0,T];H^2(\calo)),\end{equation} and the lemma is proved. $\Box$\bk

  Lemmas \ref{l5.1nou}, \ref{l5.4nou} and part (ii) now imply part (iii) and the proof of Proposition \ref{l4.2} is complete. $\Box$

\bk\n{\bf Proof of Theorem \ref{t3.1} (continued).} It is enough to prove the existence for   initial conditions $x\in H^1_0(\calo)$, provided one can also prove \eqref{e35} for such solutions with initial conditions $x,x^*\in H^1_0(\calo)$. Indeed, if we have that we can extend our solutions for arbitrary $x\in L^2(\calo)$, since $H^1_0(\calo)$ is dense in $L^2(\calo)$ and \eqref{e3.2} is obviously stable under taking limits in $X_n$ replacing $X$, with $X_n$ converging in $L^2(\ooo;C([0,T];L^2(\calo)))$ (since $\phi$ is lower semicontinuous on $L^2(\calo)$).

Hence, let $x\in H^1_0(\calo)$. Using the It\^o formula in \eqref{e5.1a} (or, equivalently, in \eqref{e5.4a}), we obtain that
$$\barr{ll}
\dd\E|X_\lbb(t)|^2_2\!\!
&\dd+2\E\int^t_0\int_\calo j_\lbb(\na(X_\lbb(s,\xi)))d\xi\,ds \dd+\lbb\E\int^t_0
|\na(X_\lbb(s))|^2_2ds\vsp
&\le|x|^2_2+\dd\frac12\ \E\dd\int^t_0\int_\calo\sum_{j=1}^\9\mu^2_j| X_\lbb e_j|^2
d\xi\,ds,\ t\in[0,T],\earr$$ because $\wt\psi_\lbb(u)\cdot u \ge j_\lbb(u) +  \lbb |u|^{2}, \forall u\in \rr^N.$

This  yields (via Gronwall's lemma)
\begin{equation}\label{e4.17}
\barr{l}
\E|X_\lbb(t)|^2_2+
2\E\dd\int^t_0\int_\calo j_\lbb(\na X_\lbb(s,\xi))d\xi\,ds\vsp
\qquad+\lbb\E  \dd\int^t_0\int_\calo|\na X_\lbb|^2 d\xi\,ds  \le C_1,\ \ \ff\lbb>0,\ t\in[0,T],\earr\end{equation}
where  $C_1=e^{2C^2_\9}|x|^2_2$.


  Moreover, we~have, for all $t\in[0,T]$,
\begin{equation}\label{e4.20}
\E \dd\int^t_0 \phi(X(t))dt\le
\dd\liminf_{\lbb\to0}\E\dd\int^t_0   \int_\calo j_\lbb(\na(X_\lbb(t))) d\xi\,dt<\9,\end{equation}
where $\phi$ is defined by \eqref{e2.4}. Indeed, we have
\begin{equation}\label{4.23}
|j_\lbb(\na u)-|\na u|\,|\le \frac12\ \lbb,\end{equation}
and this yields
\begin{equation}\label{e4.22}
\barr{r}
\dd\left|\E\left(\dd\int^t_0 \int_\calo j_\lbb(\na X_\lbb(t))d\xi\,dt-
\dd\int^t_0 \phi(X_\lbb(t))dt\right)\right|\le C\lbb,\ \ff\lbb\in(0,1].\earr\end{equation}

On the other hand, we have
\begin{equation}\label{4.25}
\lim_{\lbb\to 0}  \E\left[\sup_{t\in[0,T]} |X_\lambda(t) - X(t)|^2_2\right] = 0.\end{equation}
Indeed, by It\^o's formula,  we have
$$\barr{l}
\dd\frac 12\ d|X_\lbb(t)-X_\vp(t)|^2_2 \vsp
\qquad+ \<\psi_\lambda (\na X_\lambda(t)) -
\psi_\vp (\na   X_\vp (t)),\nabla (X_\lambda(t)- X_\epsilon(t))  \> \vsp
\qquad\dd + \<\lambda \na X_\lambda(t) -\vp \na   X_\epsilon(t) , \nabla(X_\lambda(t)-X_\vp (t)) \>  \vsp
\qquad=  \dd\frac12\int^t_0 \int _{\calo}\sum^\9_{j=1}
\mu^2_j|(X_\lbb(s)-X_\vp(s))e_j|^2ds\,d\xi\vsp
\qquad+\dd\int^t_0  \<X_\lbb-X_\vp,(X_\lbb-X_\vp)dW(s)\>,\ t\in[0,T].  \earr$$
Taking into account that, by the definition of $\psi_\lbb$,
$$ (\psi_\lambda (u) - \psi_\epsilon (v))
\cdot (u-v)  \geq (\lambda \psi_\lambda(u) - \epsilon\psi_\vp (v))\cdot (\psi_\lambda(u) -
\psi_\vp (v))  \geq -(\lambda + \vp)$$
and that
$$\barr{ll}
\<\lbb\na X_\lbb(t)-\vp\na X_\vp(t),
\na(X_\lbb(t)-X_\vp(t))\>\vsp
\qquad=-\<\lbb\D X_\lbb(t)-\vp\D X_\vp(t),
X_\lbb(t)-X_\vp(t)\>\vsp
\qquad\ge-(\lbb^2|\D X_\lbb(t)|^2_2+\vp^2|\D X_\vp(t)|^2_2)-\dd\frac12\,|X_\lbb-X_\vp|^2_2,\earr$$
we get, for some constant $C>0$ and all $t\in[0,T]$,
$$\barr{l}
  \dd |X_\lbb(t)-X_\vp(t)|^2_2\le (C^2_\9+1)\int^t_0|X_\lbb(s)-X(s)|^2_2ds+M_{\lbb,\vp}(t)\vsp
  \qquad\dd+2(\lbb+\vp)t\int_\calo d\xi+2\lbb^2\int^t_0|\D X_\lbb(s)|^2_2ds
  +2\vp^2\int^t_0|\D X_\vp(s)|^2_2 ds,\earr$$
where
$$M_{\lbb,\vp}(t)=2\int^t_0
\<X_\lbb -X_\vp ,(X_\lbb-X_\vp)dW(s)\>,\ t\in[0,T],$$
is a local real-valued $(\calf_t)$-martingale. Then, by the Burkholder-Davis-Gundy inequa\-lity (for $p=1)$, we get (see \cite{3aa}, (3.12)-(3.13)), for some constant $C>0$,
$$\barr{l}\dd\E\sup_{0\le s\le t}  |X_\lbb(s)-X_\vp(s)|^2_2
\le C(\lbb+\vp)+C\dd\int^t_0\E\sup_{0\le s\le t}| X_\lbb(s)-X_\vp(s)|^2_2ds\vsp
\qquad+C\lbb^2\E\dd\int^t_0|\D X_\lbb(s)|^2_2ds
+C\vp^2\E\dd\int^t_0|\D X_\vp(s)|^2_2ds,\
t\in[0,T],\earr$$

\noindent and, by Lemma \ref{l5.1nou} and Gronwall's lemma, it follows that $\{X_\lambda \}_\lbb$  is Cauchy in     $ L^2(\ooo;C([0,T]; $ $ L^2(\calo)))$,  which completes the proof of \eqref{4.25}.

Now, recalling that $\phi$ is lower-semicontinuous in $L^1(\calo)$ (see \eqref{e2.3}), we have by \eqref{4.25} and Fatou's lemma that
$$\liminf_{\lbb\to0}\E\int^t_0\phi(X_\lbb(t))dt\ge\E\int^t_0\phi
(X(t))dt,\ \ff t\in[0,T],$$
which, by virtue of \eqref{e4.22}, implies \eqref{e4.20}, as claimed.

We note that \eqref{4.25} and \eqref{e54prim} imply \eqref{e35}, and that \eqref{e33prim} then follows from \eqref{e54prima} and Fatou's lemma.

It remains to   prove  \eqref{e3.2}. By It\^o's formula,   we have, for all the processes $ Z$ satisfying Definition \ref{d3.1}(iii) and \eqref{e3.1}, (cf. Remark \ref{r51prim}),
\begin{equation}\label{e4.24}
\barr{l}
\dd\frac12\,\E|(X_\lbb(t)- Z(t))|^2_2+
\E\dd\int^t_0\int_\calo j_\lbb(\na X_\lbb(\tau)) d\xi\,d\tau\vsp
\qquad\le\dd\frac12\,\E|x-  Z(0)|^2_2+\E\dd\int^t_0\int_\calo j_\lbb(\na  Z(\tau)) d\xi\,d\tau\vsp
\qquad+\dd\frac12\,\E\int^t_0\<X_\lbb(\tau)-Z(\tau),G(\tau)\>d\tau\vsp
\qquad+\dd\frac12\ \E\int^t_0\int_\calo \mu(X_\lbb(\tau)-  Z(\tau))^2d\xi\,d\tau,\  t\in{[0,T]}.\earr\end{equation}

Now, letting $\lbb$ tend to zero, it follows by \eqref{e4.20}, \eqref{4.23} and \eqref{4.25}  that \eqref{e3.2} holds. This completes the proof of the existence. $\Box$

\bk\n{\bf Uniqueness.} Let $X^*$ be an arbitrary variational solution to \eqref{e3.1}   with $X^*(0)=x^*\in L^2(\calo)$ and satisfying \eqref{e34}.

 Let $x\in H^1_0(\calo)$ and $X$ be the solution constructed in the existence part of the proof, but with $X(0)=x$. Set $Y^*:=e^{-W}X^*$ and $Y:= e^{-W}X$. We set $ Y^{\vp}_\lambda= J_\vp  (Y_{\lambda}),$
where $Y_\lbb$ is the
solution to \eqref{e4.5}, but with initial condition $x\in H^1_0(\calo)$. On the basis of (H2) and Lemma \ref{l5.1nou} it follows that $e^{-W}Y^\vp_\lbb\in L^2([0,T]\times\ooo;H^1_0(\calo))$. Clearly, it is also a $\pas$ continuous $(\calf_t)$-adapted process in $L^2(\calo)$. Hence,  in \eqref{4.3a}, \eqref{e4.4}, we may choose $\wt Z=Y^\vp_\lbb$ and
we obtain that for $$G= G^{\vp}_{\lambda}= - J_\vp  ({\rm div}\ \wt\psi_\lbb(\na(e^W Y_\lbb))) +\eta^{\vp}_{\lambda},$$
where $$\barr{l}
 \eta^{\vp}_{\lambda}=\dd\frac12\,e^W(J_\vp (\mu Y_\lbb)
-\mu J_\vp (Y_{\lambda}))\vsp\qquad
 +J_\vp({\rm div}\,\wt\psi_\lbb(\na(e^WY_\lbb)))
-e^WJ_\vp(e^{-W}{\rm div}\,\wt\psi_\lbb(\na(e^WY_\lbb))),\earr$$  the function $\wt Z$ satisfies \eqref{4.3a}.\mk

 Then, by \eqref{e4.4}, we have

\begin{equation}\label{e4.25}
\barr{l}
\dd\frac12\,\E|e^{W(t)}(Y^*(t)-Y^{\vp}_\lbb(t))|^2_2
+\E\dd\int^t_0\phi(e^{W(\tau)}Y^*(\tau))d\tau\vsp
\qquad\le\dd\frac12\,|x^*-x|^2_2+\E\dd\int^t_0\phi(e^{W(\tau)}Y^{\vp}_\lbb(\tau))d\tau
\vsp
\qquad+\dd\frac12\,\E\dd\int^t_0\int_\calo\mu e^{2W(\tau)}(Y^*(\tau)-Y^{\vp}_\lbb(\tau))^2d\xi\,d\tau
\vsp
\qquad+\E\dd\int^t_0\<e^{W(\tau)}(Y^*(\tau)-Y^{\vp}_\lbb(\tau)), G^{\vp}_{\lambda}\>
d\tau,\vsp\hfill  \mbox{a.e. } t\in{[0,T]},\ \lbb>0.\earr\hspace*{-4mm}\end{equation}

 By Green's formula, we have
$$\barr{l}
\<e^W(Y^*-Y^{\vp}_{\lambda}) , G^{\vp}_\lambda \> \vsp
=\<\psi_\lambda(\na(e^{W} Y_\lambda)) + \lambda \na(e^W Y_\lambda),
\na J_\vp (e^{W} Y^*)-\na(e^WY_\lbb)\>   \vsp
+\<\psi_\lambda(\na(e^W Y_\lambda))+
\lambda \na(e^W Y_\lambda), \zeta^{\vp}_\lambda\> +
\<e^{W} (Y^*-Y^{\vp }_\lambda),\eta^{\vp }_\lambda\>\earr$$ where
$$ \zeta_\lambda^{\vp } =    \na(e^{W}Y_{\lambda})- \na  J_\vp (e^{W}Y^{\vp }_\lambda ) . $$
 Taking into account that
$$\psi_\lbb(u)\cdot(u-v)\ge j_\lbb(u)-j_\lbb(v),\ \ff u,v\in \rr^d,$$
this yields
$$\barr{l} \<e^W (Y^*- Y^{\vp }_\lambda), G^{\vp }_\lambda \>
 \leq    \phi_\lambda
 (J_\vp (e^{W} Y^*) )
 -\phi_\lbb(e^WY_\lbb)-
\lbb|\na(e^WY_\lbb)|^2_2 \vsp
\qquad-\lambda\<\D(e^W Y_\lambda),
J_\vp (e^{W} Y^*)\>  + \<e^W(Y^*-Y^{\vp }_\lambda),\eta^{\vp }_\lambda\>\vsp
\qquad+
 \< \psi_\lbb(\na(e^W Y_\lambda))+ \lambda \na(e^W Y_\lambda),
 \zeta^{\vp }_\lambda\>.\earr$$
Here,   $\phi_\lbb$ is the function
$$ \phi_\lbb(z)=\int_\calo j_\lbb(\na z)d\xi,\ \ \ff z\in H^1_0(\calo) .$$
 Substituting into \eqref{e4.25}, we obtain that

\begin{equation}\label{e4.26}
\barr{l}\dd\frac12\,\E|e^{W(t)}(Y^*(t)-Y^{\vp }_\lbb(t))|^2_2+
\E\dd\int^t_0\phi (e^{W(\tau)}Y^*(\tau))d\tau\\
\quad +\E\dd\int^t_0\phi_\lbb(e^{W(\tau)}
Y_\lbb(\tau))d\tau+\lbb\E\dd\int^t_0
|\na(e^{W(\tau)}Y_\lbb(\tau))|^2_2d\tau\\
 \quad\le\dd\frac12\,|x^*-x|^2_2+
 \E\dd\int^t_0\phi(e^{W(\tau)}Y^\vp_\lbb(\tau))d\tau\\
 \quad\dd
+\E\dd\int^t_0\phi_\lbb ( J_\vp (e^{W}Y^*(\tau)))d\tau\\
\quad -\lambda \E\dd\int^t_0\<\D(e^{W(\tau)}Y_\lbb(\tau)),
J_\vp (e^{W(\tau)} Y^*(\tau))\>d\tau\vsp
\quad\dd+\E\int^t_0\(\<e^W(Y^*- Y^{\vp }_\lambda),\eta^{\vp }_\lbb\> + \<\psi_\lambda(\na(e^W Y_\lambda)) + \lambda \na(e^W Y_\lambda), \zeta^{\vp }_\lambda\>\)d\tau\vsp
\quad +\dd\frac12\,C^2_\9\E\dd\int^t_0(e^{W(\tau)}
(Y^*(\tau)-Y^\vp_\lbb(\tau)))|^2_2d\tau,\ t\in[0,T],\ \ff\lbb>0,\earr\hspace*{-15mm}\end{equation}
 where  $C^2_\9$ is as in (H1).
Now, as seen earlier in \eqref{e4.22}, we have
\begin{equation}\label{e4.27}
|\phi(e^{W(\tau)}Y_\lbb(\tau))-\phi_\lbb(e^{W(\tau)}Y_\lbb(\tau))|\le C\lbb,\ \ \ff\tau\in[0,T].\end{equation}
Similarly,  we have also
 \begin{equation}\label{e4.28}
\dd\int^T_0\!\!|\phi_\lambda(J_\vp (e^{W(\tau)}Y^*(\tau))) - \phi ( J_\vp (e^{W(\tau)}Y^*(\tau))|d\tau
\dd\leq  C\lambda .\end{equation}
Substituting \eqref{e4.27}, \eqref{e4.28} in \eqref{e4.26},   yields
 \begin{equation}\label{e4.33a}
 \barr{l}
 \dd\frac12\,\dd \E |e^{W(t)}(Y^*(t) -Y^{\vp }_\lambda(t))|^2_2  +\E\dd\int^t_0\phi(e^{W(\tau)}Y^*(\tau))d\tau\vsp
  \quad\leq \dd\frac12\,|x^*-x|^2_2+
  \E\dd\int^t_0\phi(J_\vp(e^{W(\tau)}Y^*(\tau)))d\tau\vsp
 \quad\dd
  +\E\dd\int^t_0(\phi(e^{W(\tau)}Y^\vp_\lbb(\tau))
  -\phi(e^{W(\tau)}Y_\lbb(\tau)))d\tau\vsp
  \quad+\dd\frac12\,C^2_\9
  \E\dd\int^t_0|e^{W(\tau)}(Y^*(\tau)-Y^\vp_\lbb(\tau))|^2_2d\tau\vsp
  \quad-\lbb\E\dd\int^t_0\<\D(e^{W(\tau)}Y_\lbb(\tau)),
  J_\vp(e^{W(\tau)}Y^*(\tau))\>d\tau\vsp
  \quad+C_{\lbb,\vp}
  \(\E\dd\int^t_0    |\zeta^\vp_\lbb(\tau)|^2_2 d\tau\)^{1/2}
  +C_{\lbb,\vp}
  \(\E\(\dd\int^t_0|\eta^\vp_\lbb(\tau)|^2_2d\tau\)^{r/2}\)^{1/r},
 \earr\end{equation}
 where $\de$ is as in \eqref{e34}, $r=\frac{\de+2}{\de+1}$ and
 $$\barr{l}
 \dd C_{\lbb,\vp}=4\(\E\(\dd\int^T_0
 |e^W(Y^*-Y^\vp_\lbb)|^2_2d\tau\)^{\frac{2+\de}2}\)^{\frac1{2+\de}}
 \vsp
 \qquad\dd+4+
 4\(\E\dd\int^T_0 \lbb|\na(e^WY_\lbb)|^2_2d\tau\)^{1/2}.\earr$$
 Now, recalling that, by Corollary \ref{c9.1},
 $$\E\int^t_0\phi(J_\vp(e^{W(\tau)}Y^*(\tau)))d\tau
 \le\E\int^t_0\phi(e^{W(\tau)}Y^*(\tau))d\tau,\ \ff\vp>0,$$
 letting $\vp\to0$ in  \eqref{e4.33a} yields
 \begin{equation}\label{5.35}
 \barr{l}
 \E|e^{W(t)}(Y^*(t)-Y_\lbb(t))|^2_2\le |x^*-x|^2_2
 \vsp\qquad\dd+\,
 C^2_\9\E\dd\int^t_0|e^{W(\tau)}
 (Y^*(\tau)-Y_\lbb(\tau))|^2d\tau\vsp
\qquad\dd -\,\lbb\E\dd\int^t_0\<\D(e^{W(\tau)}Y_\lbb(\tau)),
 e^{W(\tau)}Y^*(\tau)\>d\tau.\earr\end{equation}
 because

 $$\barr{l}
 \dd\lim_{\vp\to0}\E
\(\dd\int^T_0
 |\eta^\vp_\lbb(\tau)|^2_2d\tau\)^{r/2}=0,\vsp
\dd\lim_{\vp\to0}\E\dd\int^T_0
|\zeta^\vp_\lbb(\tau)|^2_2 d\tau=0,\vsp
 \sup\{C_{\lbb,\vp};\ \vp\in(0,1)\}<\9,\vsp
 e^WY_\lbb=X_\lbb\in L^2([0,T]\times\ooo;H^2(\calo))\mbox{ by Lemma \ref{l5.1nou} and}\vsp
 e^WY^*\in L^2([0,T]\times\ooo;L^2(\calo)).\earr$$
To check all this is pretty routine. The main problem is to justify the interchange of "$\lim_{\vp\to0}$" with the integral with respect to $d\tau\otimes\mathbb{P}$, i.e., to find an integrable uniformly dominating function. As an exemplary case, we show how this is done for the  last summand in the definition of $\eta^\vp_\lbb$:

Clearly, since $J_\vp$ is a contraction on $L^2(\calo)$, it follows by \eqref{e521prim} that there exists a constant $c>0$ such that, for all $\vp\in(0,1)$,
$$|e^WJ_\vp(e^{-W}{\rm div}\,\wt\psi_\lbb(\na(e^WY_\lbb)))|^2_2\le c\cdot\exp\(4\sup_{\tau\in[0,T]}|W(\tau)|_\9\)
\|X_\lbb\|^2_{H^2(\calo)}.$$Hence, applying H\"older's inequality with $p=\frac2r\ (>1),$ $q=\frac2{2-r}$ to the expectation, we obtain
$$\barr{l}
\E\(\dd\int^T_0\exp\(4\sup_{\tau\in[0,T]}|W(\tau)|_\9\)
\|X_\lbb\|^2_{H^2(\calo)}d\tau\)^{r/2}\vsp
\qquad\dd\le\(\E\exp\(\frac{8r}{2-r}\sup_{\tau\in[0,T]}
|W(\tau)|_\9\)\)^{\frac{2-r}2}
\(\E\dd\int^T_0\|X_\lbb(\tau)\|^2_{H^2(\calo)}d\tau\)^{r/2},\earr$$which is finite by \eqref{e1.2prim} and Lemma \ref{l5.1nou}.

 Now, by Lemma \ref{l5.1nou}, we have
 $$\lim_{\lbb\to0}\lbb\E\int^t_0\<\D(e^{W(\tau)}Y_\lbb(\tau)),
 e^{W(\tau)}Y^*(\tau)\>d\tau=0.$$
 Then, letting $\lbb\to0$ in \eqref{5.35}, we obtain via Gronwall's lemma
 $$\E|X^*(t)-X(t)|^2_2=
 \E|e^{W(t)}(Y^*(t)-Y(t))|^2_2
 \le|x^*-x|^2_2 e^{C^2_\9T}.$$
 Now, letting $x\to x^*$ in $L^2(\calo)$, we see by \eqref{e35}  that $X^*$ coincides with the solution starting at $x^*$ constructed in the existence part of the proof, which is hence unique. $\Box$

\begin{remark}\label{r5.6} {\rm We did not succeed in proving the uniqueness for Theo\-rem~\ref{t3.1} directly for the original equation \eqref{e1.1}. The reason is that, regularizing \eqref{e1.1} by $J_\vp$ destroys the special form of the noise. Therefore, we had to use equation \eqref{e4.1} and Proposition \ref{p4.2}.}\end{remark}

\section{Positivity of solutions}
\setcounter{equation}{0}
\setcounter{theorem}{0}

It should be emphasized that physical models of nonlinear  diffusion are concerned in general with nonnegative solutions of equation \eqref{e1.1}. In this context, we have  the following result.

\begin{theorem}\label{t5.1} In Theorem {\rm\ref{t3.1}} assume further that $x\ge0$, a.e. in $\calo$. Then
\begin{equation}\label{e5.1}
X(t,\xi)\ge0\mbox{\ \ a.e. in }(0,T)\times\calo\times\ooo.\end{equation}\end{theorem}

\n{\bf Proof.}  It suffices to show that the solution $X_\lbb$ to \eqref{e5.1a} is a.e. nonnegative on $[0,T]\times\calo\times\ooo$. By \eqref{e54prim} we may assume that $x\in L^4(\calo)$. Below we only give a heuristic argument to prove the assertion (e.g., apply It\^o's formula in an informal way), which can be made rigorous by regularization. Since the latter is analogous as in the proof of Theorem 2.2 in \cite{5prim} or can be done similarly as in the proof of Theorem \ref{t6.1} below, we omit the details.

We apply the It\^o formula in \eqref{e5.1a} to the function $x\to\frac14\ |x^-|^4_4.$ We obtain
$$\barr{l}
\dd\frac14\ \E\int_\calo|X_\lbb^-(t,\xi)|^4d\xi+\E\int^t_0
\int_\calo\wt\psi_\lbb(\na X_\lbb(s,\xi))\cdot\na|X_\lbb^-(s,\xi)|^3d\xi\,ds\vsp
\qquad=\dd\frac14\int_\calo|x^-(\xi)|^4d\xi+\E\dd\int^t_0\int_\calo\dd\sum^\9_{j=1}\mu^2_j(X_\lbb e_j)^2(X_\lbb^-)^2d\xi\,ds.\earr$$

Recalling that $\na y\cdot\na y^-=-|\na y^-|^2$ a.e. in $\calo$ for each $y\in H^1(\calo)$, it follows that
$$\barr{r}
\E\dd\int_\calo |X^-_\lbb(t,\xi)|^4d\xi\le
 C\E\dd\int^t_0\int_\calo |X^-_\lbb(t,\xi)|^4d\xi\,ds, \ff t\in[0,T],\earr$$
which implies that $X^-_\lbb\equiv0$, as claimed.

\section{Extinction in finite-time}
\setcounter{equation}{0}
\setcounter{theorem}{0}

A striking feature of highly singular nonlinear diffusion equations is the extinction in finite time of the solution. In nonlinear diffusion phenomena, this is due to the   singularity at level $X=0$ of the diffusivity
and this causes a fast loss of mass. (See \cite{3a} for the case of stochastic porous media equation and \cite{3aa}, \cite{6a} for stochastic self-organized criticality.)
A similar phenomenon happens in the case of equation \eqref{e1.1}.

\begin{theorem}\label{t6.1} Let $2 \le  N  \le   3.$  Let $X$ be as in Theorem {\rm \ref{t3.1}}, with initial condition $x\in L^N(\calo)$, and let \mbox{$\tau{=}\inf\{t\ge0; |X(t)|_N=0\}$.} Then, we have
\begin{equation}\label{e6.1}
\mathbb{P}[\tau\le t]\ge 1-\rho\1\(\int^t_0 e^{-C^*s} ds\)^{-1}|x|_N,\ \ff t\ge 0. \end{equation}
Here $\rho=\inf\{|y|_{W^{1,1}_0(\calo)}/|y|_{\frac N{N-1}}$; $ y\in W^{1,1}_0(\calo)\}$ and $C^*= \frac{C^2_\9}2\ N(N-1).$
In particular, if $|x|_N<\rho/C^*$, then $\mathbb{P}[\tau<\9]>0.$
\end{theorem}

We shall prove Theorem \ref{t6.1} as stated, i.e., only for $2\le N\le3.$ The case $N=1$ is similar, but one proves extinction in $L^2(\calo)$-norm rather than $L^1(\calo)$-norm (see \cite[Theorem 3]{6aa} for details). We fix $\lbb\in(0,1]$ and start with the \fwg\ lemma, which is one of the main ingredients of the proof.

Before, we recall that, by \eqref{5.5a}, $X_\lbb$ is $\pas$ continuous in $H^1_0(\calo)$. For $K\in\nn,$ $K>\|x\|_1$, define the $\{\calf_t\}$-stopping time
$$\tau_K:=\inf\{t\ge0; \|X_\lbb(t)\|_1>K\}.$$

\begin{lemma}\label{l7.2} Let $x\in H^1_0(\calo)$. For every $K\in\nn$, $K>\|x\|_1$, we have $\pas$
\begin{equation}\label{ez7.2}
\barr{l}
|X_\lbb(t)|^N_N+N\rho
\dd\int^t_s|X_\lbb(r)|^{N-1}_Ndr\vsp
\le|X_\lbb(s)|^N_N+C^*
\dd\int^{t}_{s}|X_\lbb(r)|^N_Ndr
+N(N-1)\lbb\int^t_s|X_\lbb(r)|^{N-2}_{N-2}dr\vsp
+N\dd\int^{t}_{s}
\<|X_\lbb(r)|^{N-2}X_\lbb(r),
X_\lbb(r)dW(r)\>,\ \ff s,t\in[0,T],\ s\le t.\earr\end{equation}
\end{lemma}

\n{\bf Proof.} Since $N\le3,$ we have by Sobolev embedding, $H^1_0(\calo)\subset L^4(\calo)$ continuously, hence, for some constant $C>0$,
\begin{equation}\label{ez7.3}
\sup_{t\in[0,\tau_K]}|X_\lbb(t)|_N\le CK\mbox{ on }\ooo.\end{equation}

We have by standard interpolation (see, e.g., \cite[Theorem 2.1]{29})  if $N=3$
$$\barr{ll}
\dd\E\int^{\tau_K}_0\|X_\lbb(t)\|^3_{1,3}dt\!\!\!
&\le C\E\dd\int^{\tau_K}_0\|X_\lbb(t)
\|^2_{H^2(\calo)}|X_\lbb(t)|_3dt\vsp
&\le CK\E\dd\int^T_0\|X_\lbb(t)\|^2_{H^2(\calo)}dt<\9\earr$$by \eqref{e57prim},  and  if $N=2$
$$\E\int^{\tau_K}_0\|X_\lbb(t)\|^2_1dt<\9.$$ Hence, by Theorem 2.1 in \cite{25prim}, applied with
$$\barr{rcl}
f_t&:=&\wt\psi_\lbb(\na X_\lbb(t))\ (\le1+\lbb|\na X_\lbb(t)|)\vsp
f^\circ_t&:=&0\vsp
g^h_t&:=&\mu_ke_kX_\lbb(t),\earr$$we have the following It\^o formula for the $L^N(\calo)$-norm $\pas$
\begin{equation}\label{ez7.4}\!\!\!\barr{l}
|X_\lbb(t\wedge\tau_K)|^N_N+N(N-1)
\dd\int^{t\wedge\tau_K}_{s\wedge\tau_K}\int_\calo|X_\lbb(r)|^{N-2}
\na X_\lbb(r)\cdot\wt\psi_\lbb(\na X_\lbb(r))d\xi\,dr\vsp
\quad=|X_\lbb(s\wedge\tau_K)|^N_N+\dd\frac12\,N(N-1)
\dd\int^{t\wedge\tau_K}_{s\wedge\tau_K}\int_\calo\mu|X_\lbb(r)|^Nd\xi\,dr\vsp
\quad\dd+N\dd\int^{t\wedge\tau_K}_{s\wedge\tau_K} \<|X_\lbb(r)|^{N-2}X_\lbb(r),X_\lbb(r)dW(r)\>,\ \ff s,t\in[0,T],\ s\le t.\earr\hspace*{-15mm}\end{equation}


Since, by interpolation (cf. \cite[Theorem 2.1]{29})
$$\barr{lcl}
\dd\E\int^T_0|X_\lbb(r)|^6_3dr&\le&\dd
C\E\int^T_0(\|X_\lbb(r)\|^{\frac32}_{H^2(\calo)}
|X_\lbb(r)|^{\frac92}_2)dr\vsp
&\le&C\E\dd\int^T_0(\|X_\lbb(r)\|^2_{H^2(\calo)}+
|X_\lbb(r)|^{9}_2)dr,\earr$$and, since by \eqref{e57prim} and \eqref{e54prima} the last term is finite,   we can let $K\to\9$ in \eqref{ez7.4} to obtain

\begin{equation}\label{ez74prim}
\barr{l}
|X_\lbb(t)|^N_N+N(N-1)\dd\int^t_s\int_\calo
|X_\lbb(r)|^{N-2}\na X_\lbb(r)\cdot\wt\psi_\lbb(\na X_\lbb(r))d\xi\,dr\vsp
\qquad=\dd|X_\lbb(s)|^N_N+\frac12\,N(N-1)
\dd\int^t_s\int_\calo\mu|X_\lbb(r)|^Nd\xi\,dr\vsp
\qquad+N\dd\int^t_s
\<|X_\lbb(r)|^{N-2}
X_\lbb(r),X_\lbb(r)dW(r)\>,\ \ff s,t\in[0,T],\ s\le t.\earr\hspace*{-15mm}\end{equation}

But, since $\wt\psi_\lbb(u)\cdot u\ge|u|-\lbb,$ we have
\begin{equation}\label{ez7.5}\barr{l}
(N-1)|X_\lbb|^{N-2}\na X_\lbb(r)\cdot\wt\psi_\lbb(\na X_\lbb)\vsp\qquad
\ge(N-1)|X_\lbb|^{N-2}(|\na X_\lbb|-\lbb)\vsp\qquad
=|\na(|X_\lbb|^{N-1})|
-(N-1)\lbb|X_\lbb|^{N-2}.\earr
\end{equation}
Hence, the second term on the left hand side of \eqref{ez7.4} is bigger than
$$N\rho
\dd\int^{t}_{s}\int_\calo|X_\lbb(r)|^{N-1}_Ndr
-N(N-1)\lbb\int^t_s|X_\lbb(r)|^{N-2}_{N-2}dr,$$
where we used Sobolev's embedding theorem in $W^{1,1}_0(\calo)$, i.e.,
$$\rho|y|_{\frac N{N-1}}\le\|y\|_{1,1},\ \ff y\in W^{1,1}_0(\calo),$$in the last step. Plugging this into \eqref{ez7.4}  implies the assertion of the lemma.~$\Box$

\begin{lemma}\label{l7.3} Let $x,y\in H^1_0(\calo)$ and let $X^x_\lbb,X^y_\lbb$ denote the solutions to \eqref{e5.1a} with initial conditions $x,y$, respectively. Then $\pas$
\begin{equation}\label{ez75prim}
 \!\!\barr{ll}
|X^x_\lbb(t){-}X^y_\lbb(t)|^N_N \dd\le
|X^x_\lbb(s){-}X^y_\lbb(s)|^N_N
+C^*\int ^t_s
|X^x_\lbb(r){-}X^y_\lbb(r)|^N_Ndr\vsp
 \dd{+}N\!\int^t_s\!\<|X^x_\lbb(r){-}X^y_\lbb(r)|^{N-2}
(X^x_\lbb(r){-}X^y_\lbb(r)),
(X^x_\lbb(r){-}X^y_\lbb(r))dW(r)\>,\vsp
 \hfill\ff s,t\in[0,T],\ s\le t.\earr\end{equation}
Furthermore, for some constant $C$ independent of $\lbb$
\begin{equation}\label{ez75secund}
\E\left[\sup_{t\in[0,T]}|X^x_\lbb(t)-X^y_\lbb(t)|^N_N\right]
\le2|x-y|^N_Ne^{CT},\ \ff x,y\in H^1_0(\calo).\end{equation}
\end{lemma}

\n{\bf Proof.} \eqref{ez75prim} follows analogously to \eqref{ez74prim}, taking into account that $(\wt\vf_\lbb(u)-\wt\vf_\lbb(v))\cdot(u-v)\ge0$, $\ff u,v\in\rr^N$. Then \eqref{ez75secund} follows by a standard application of the Burkholder-Davis-Gundy inequality (for $p=1)$. $\Box$\bk


\n{\bf Proof of Theorem \ref{t6.1} (continued).} Let $x\in H^1_0(\calo)$ and $X^x_\lbb$ be the solution to \eqref{e5.1a} with initial condition $x$. Applying It\^o's formula to \eqref{ez74prim}  and the function $\vf_\vp(r)=(r+\vp)^{\frac1N},$ $\vp\in(0,1)$, and proceeding as in the proof  of the previous lemma, we obtain $\pas$
\begin{equation}\label{ez75tert}
\barr{l}
\vf_\vp(|X^x_\lbb(t)|^N_N)+\rho
\dd\int^{t}_{s}|X^x_\lbb(r)|^{N-1}_N
(|X^x_\lbb(r)|^N_N+\vp)^{-\frac{N-1}N}dr\vsp
\qquad\le\vf_\vp(|X^x_\lbb(s)|^N_N)+C^*
\dd\int^{t}_{s}|X^x_\lbb(r)|_Ndr\vsp
\qquad+\lbb(N-1)\dd\int^t_s|X^x_\lbb(r)|^{N-2}_{N-2}
(|X^x_\lbb(r)|^N_N+\vp)^{-\frac{N-1}N}dr\vsp
\qquad+\dd\int^{t}_{s}\<X^x_\lbb(r)|X^x_\lbb(r)|^{N-2}
(|X^x_\lbb(r)|^N_N+\vp)^{-\frac{N-1}N},X^x_\lbb(r)dW(r)\>,\vsp
\hfill\ff s,t\in[0,T],\ s\le t.\earr\end{equation}
Hence, one can let $\vp\to0$ in \eqref{ez75tert} to arrive at
\begin{equation}\label{e10}
\!\!\barr{l}
e^{-C^*t}|X^x_\lbb (t)|_N+\rho\dd\int^t_s   \one_{[|X^x_\lbb(\theta)|_N>0]}e^{-C^*\theta} d\theta  \dd\le
e^{-C^*s}|X^x_\lbb(s)|_N\vsp\dd+\dd\int^t_s    \one_{[|X^x_\lbb(\theta)|_N>0]}e^{-C^*\theta} |X^x_\lbb(r)|^{-(N-1)}_N
\<X^x_\lbb(r)|X^x_\lbb(r)|^{N-2},X^x_\lbb(r)dW(r)\>.
\earr\hspace*{-10mm}\end{equation}
In particular, this implies that the process $t\to e^{-C^{*}t}|X^x_\lbb(t)|_N$ is an $\{\calf_t\}$-super\-mar\-tin\-gale.

 If in \eqref{e10} we take expectation and set $s=0$, we see that
\begin{equation}\label{e711}
e^{-C^{*}t}\E|X^x_\lbb(t)|_N
+\rho\int^t_0e^{-C^{*}\theta}
\PP[|X^x_\lbb(\theta)|_N>0]d\theta\le|x|_N,\ \ff t>0.\end{equation}
Since $x\in H^1_0(\calo)$, by Lemma \ref{l5.1nou}, Remark \ref{r3.6} and interpolation we have,  for $N=3$ and some $C>0$,
\begin{equation}\label{e712}
\barr{l}
\dd\E\left[\sup_{t\in[0,T]}
|X^x_\lbb(t)-X^x(t)|^2_N\right]\vsp\qquad
\le\dd C\(\E\left[\sup_{t\in[0,T]}
|X^x_\lbb(t)-X^x(t)|^2_2\right]\right)^{\frac12}
\|x\|_1,\  \ff\lbb\in(0,1],\earr\end{equation}
where $X^x$ is the solution to \eqref{e1.1} with initial condition $x$.

Hence, by \eqref{4.25}, for $2\le N\le3$,
\begin{equation}\label{e713}
\lim_{\lbb\to0}\E\left[\sup_{t\in[0,T]}
|X^x_\lbb(t)-X^x(t)|^2_N\right]=0,\end{equation}

Noting that, since for each $t>0$,
$$\int^t_0e^{-C^*\theta}
\mathbb{P}[|X^x_\lbb(\theta)|_N>0]d\theta
=\dd\sup_{\vp>0}\dd\int^t_0 e^{-C^*}
\E[|X^x_\lbb(\theta)|_N(|X^x_\lbb(\theta)|_N
+\vp)^{-1}]d\theta,$$
by \eqref{e713} and Fatou's lemma, \eqref{e711} also holds with $X^x$ replacing $X^x_\lbb$.

 But then, by \eqref{e35} and Fatou's lemma, \eqref{e711} extends to every solution $X^x$ of \eqref{e1.1} for arbitrary $x\in L^2(\calo)$.
 But, for $x\in L^2(\calo)$, by \eqref{e713} and \eqref{e35}, the process $t\to e^{-C^*t}|X^x(t)|_N$ is an $L^1$-limit of supermartingales, hence itself is a supermartingale. Hence
 $$|X^x(t)|_N=0\mbox{ for }t\ge\tau=\inf\{t\ge0:
 |X^x(t)|_N=0\},$$ and thus $\mathbb{P}[|X^x(\theta)|_N>0]=\mathbb{P}[\tau>\theta]$. By \eqref{e711}, for $X^x$ with $x\in L^2(\calo)$, this yields $\mathbb{P}[\tau>t]\le\(\rho\dd\int^t_0
e^{-C^{*}\theta}d\theta\)\1|x|_N,$ as claimed. $\Box$

\begin{remark}\label{r6.1} {\rm  In particular, taking   $\mu_k=0$ for all $k,$ implying $C^*=0$, we have $\tau \le {\frac{|x|_{N}}\rho }$ and  recover the deterministic case (\cite{1aa}).}\end{remark}

\section{Appendix 1}
\setcounter{equation}{0}

 Proposition \ref{p9.1} below is due to H.~Brezis (\cite{15aa}) who answered a question we raised  and we are grateful to him for this.

\begin{proposition}\label{p9.1} Let $\calo$ be a bounded, convex domain of $\rr^N$, $N\ge1$, with smooth boundary $($of class $C^2)$. Let $J_\vp=(I+\vp A)^{-1}$, where $A=-\D$, $D(A)=H^1_0(\calo)\cap H^2(\calo)$ and $\vp>0$. Then
\begin{equation}\label{ee9.1}
\int_\calo|\na J_\vp(y)|d\xi\le\int_\calo|\na y|d\xi,\ \ \ff y\in W^{1,1}_0(\calo). \end{equation}\end{proposition}

\n{\bf Proof.} For simplicity, we shall write here $|\na y|$ instead of $|\na y|_N$. Rescaling, we can assume $\vp=1$ and so reduce \eqref{ee9.1} to
\begin{equation}\label{ee9.2}
\int_\calo|\na u|d\xi\le\int_\calo|\na y|d\xi,\ \ \ff y\in W^{1,1}_0(\calo),\end{equation}where
\begin{equation}\label{ee9.3}
u-\D u=y\mbox{ in }\calo;\ \ u=0\ dS-\mbox{a.e.  on }\pp\calo,\end{equation}and  $dS$ is the surface measure on $\pp\calo$.
 Without loss of generality, we may also assume  $y\in C^\9_0(\calo)$.

We set
\begin{equation}\label{ee9.4}
\barr{c}
\dd D_i=\frac\pp{\pp\xi_i},\ D^2_{ij}=\frac{\pp^2}{\pp\xi_i\pp\xi_j} ,\ i,j=1,...,N,\vsp
\dd\vf(\xi)=|\na u(\xi)|=\(\sum^N_{i=1}|D_iu|^2\)^{\frac12},\ \vf_\vp(\xi)=\sqrt{\vp^2+|\na u(\xi)|^2}.\earr\end{equation}

We shall prove \eqref{ee9.2} following several steps.

\begin{lemma}\label{l9.1} We have
\begin{equation}\label{ee9.5}
\frac{\vf^2}{\vf_\vp}-\D\vf_\vp\le|\na y|\ \mbox{ in }\calo.\end{equation}\end{lemma}

\n{\bf Proof.} By \eqref{ee9.4}, we have
$$\vf_\vp D_j\vf_\vp=\sum^N_{i=1}D_iuD^2_{ij}u,$$which yields
$$(D_j\vf_\vp)^2\le\frac1{\vf^2_\vp}\ |\na u|^2\sum^N_{i=1}(D^2_{ij}u)^2$$and therefore
\begin{equation}\label{ee9.6}
|\na\vf_\vp|^2\le\frac{\vf^2}{\vf^2_\vp}\ \sum^N_{i,j=1}|D^2_{ij}u|^2\le\dd\sum^N_{i,j=1}|D^2_{ij}u|^2\mbox{\ \ in }\calo.\end{equation}
We also have
$$\barr{lcl}
\vf_\vp\D\vf_\vp+|\na\vf_\vp|^2&=&
\dd\sum^N_{i,j=1}|D^2_{ij}u|^2+\sum^N_{i=1}D_iu\D D_iu\vsp
&=&\dd\sum^N_{i,j=1}|D^2_{ij}u|^2+\sum^N_{i=1}D_iu(D_iu-D_iy)\vsp
&=&\dd\sum^N_{i,j=1}|D^2_{ij}u|^2+|\na u|^2-\na u\cdot\na y\ge|\na\vf_\vp|^2+\vf^2-\vf|\na y|,\earr$$where the last  inequality follows by \eqref{ee9.6}.
This yields
$$-\vf_\vp\D\vf_\vp+\vf^2\le\vf|\na y|,$$which implies \eqref{ee9.5}, as claimed.

Assume that $0\in\pp\calo$ and represent locally
 $\pp\calo=\{(\xi',\xi_N);\ \xi_N=\gamma(\xi')\},$ where $\gamma$ is a $C^2$-function in a neighborhood of $0$ in $\rr^{N-1}$ and $\gamma(0)=0$, $\na\gamma(0)=0$.

\begin{lemma}\label{l9.2} We have
\begin{equation}\label{ee9.7}
D_N\vf_\vp(0)=(D_Nu)^2(0)(\vp^2+(D_Nu)^2(0))^{-\frac12}\D_{\xi'}\gamma(0).
\end{equation}\end{lemma}

\n{\bf Proof.} By \eqref{ee9.4}, we have
 \begin{equation}\label{ee9.8}
\vf_\vp D_N\vf_\vp=\sum^N_{i=1}D_i uD_{Ni}u. \end{equation}
Since $u=0$ on $\pp\calo$, we have
$$u(\xi_1,\xi_2,...,\xi_{N-1},\gamma(\xi_1,\xi_2,...,\xi_{N-1}))=0$$and differentiating with respect to $\xi_i$, $i=1,...,N-1$, yields
 \begin{eqnarray}
 &&D_iu+D_NuD_i\gamma\equiv0,\ i=1,...,N-1,\label{ee9.9}\\[2mm]
 &&D^2_{ii}u+2D_{iN}uD_i\gamma+D^2_{NN}u(D_i\gamma)^2+D_NuD_{ii}
 \gamma\equiv0.\label{ee9.10}
 \end{eqnarray}
  By \eqref{ee9.9}, \eqref{ee9.10}, we get in $\xi'=0,\ \xi_N=0,$
\begin{equation}\label{ee9.11}
 D_iu(0)=0,\ \
 D_{ii}u(0)+D_Nu(0)D_{ii}\gamma(0)=0.
\end{equation}
By \eqref{ee9.11}, we have
$$\D u(0)=D_{NN}u(0)-D_Nu(0)\D_{\xi'}\gamma(0),$$while, by \eqref{ee9.3}, we have
 $\D u(0)=0,$ which yields
\begin{equation}\label{ee9.13}
D_{NN}u(0)=D_N(0)\D_{\xi'}\gamma(0).\end{equation}
Now, taking \eqref{ee9.8} in $0$ and using \eqref{ee9.11}, \eqref{ee9.13}, we obtain
$$\sqrt{\vp^2+(D_Nu(0))^2}\ D_N\vf_\vp(0)=D_Nu(0)D_{NN}u(0)=(D_Nu)^2(0)\D_{\xi'}\gamma(0),$$as claimed.\bk

\n{\bf Proof of Proposition \ref{p9.1}.} Let $n$ be the outward normed to $\pp\calo$. We have
\begin{equation}\label{e912a}
\frac{\pp\vf_\vp}{\pp n}\ (0)=-D_N\vf_\vp(0)=\frac{-(D_N u(0))^2}{\sqrt{\vp^2+(D_Nu)^2(0)}}\ \D_{\xi'}\gamma(0).\end{equation}
On the other hand, since $\calo$ is convex, we have $\D\gamma(0)\ge0$ and, therefore,
\begin{equation}\label{813}
\frac{\pp\vf_\vp}{\pp n}\ (0)\le0.\end{equation}
Since $0$ can be replaced by an arbitrary point of $\pp\calo$, we have therefore
\begin{equation}\label{ee9.14}
\frac{\pp\vf_\vp}{\pp n}\le0\ \ \mbox{on }\pp\calo.\end{equation}
Integrating \eqref{ee9.5} over $\calo$, we get
$$\int_\calo\ \frac{\vf^2}{\vf_\vp} \ d\xi-\int_{\pp\calo}\ \frac{\pp\vf_\vp}{\pp n}\ dS\le\int_\calo|\na y|d\xi,$$and so, by \eqref{ee9.14}, we have
$$\int_\calo\ \frac{\vf^2}{\vf_\vp}\ d\xi\le\int_\calo |\na y|d\xi.$$Then, letting $\vp\to0$, we get \eqref{ee9.2}, thereby completing the proof.

 \begin{remark}\label{r9.1} {\rm Proposition \ref{p9.1}, which has an interest in itself, amounts to saying that the heat flow on convex smooth domains $\calo$ is nonexpansive in $W^{1,1}_0(\calo)$. Analyzing the previous proof, one sees that it remains true for domains with piecewise smooth and convex boundary.}\end{remark}

\begin{corollary}\label{c9.1} Let $\calo$ be a convex, bounded and open subset of $\rr^N$. Then
\begin{equation}\label{ee5.27}
\phi(J_\vp(u))\le\phi(u),\ \ \ff u\in BV(\calo),\end{equation}
where $\phi$ is the functional \eqref{e2.3}.\end{corollary}

\n{\bf Proof.} Let $u\in BV(\calo),\,\phi(u)<\9$. This means   that there is $\{u_n\}{\subset} W^{1,1}_0(\calo)$ such that $u_n\to u$ in $L^1(\calo)$ and
\begin{equation}\label{ee5.28}
\phi(u)\ge\limsup_{n\to\9}\int_\calo|\na u_n|d\xi.\end{equation}
while
\begin{equation}\label{ee5.29}
\limsup_{n\to\9}\int_\calo|\na J_\vp(y_n)|d\xi\ge\phi(J_\vp(y)). \end{equation}

By  \eqref{ee5.27}, \eqref{ee5.28} and \eqref{ee9.1}, it follows \eqref{ee5.27}. $\Box$\bk

Proposition \ref{p9.1} can be extended as follows.

\begin{proposition}\label{p8.6} Let $g:[0,\9)\to[0,\9)$ be a \ct\ and convex function of at most quadratic growth such that $g(0)=0$. Then
\begin{equation}\label{e8.18} \int_\calo g(|\na J_\vp(y)|)d\xi\le\int_\calo g(|\na y|)d\xi,\ \ff y\in H^1_0(\calo).\end{equation}\end{proposition}

\n{\bf Proof.} Since $g$ is of at  most quadratic growth, as before  we may assume that $y\in C^\9_0(\calo)$. Furthermore, without   loss of generality, we may assume that $g\in C^2([0,\9))$. (This can be achieved by regularizing the function $g$.) As in the previous case, it suffices to prove \eqref{e8.18} for $\vp=1$. We set
$$\phi(\xi)=g(\vf(\xi)),\ \phi_\vp(\xi)=g(\vf_\vp(\xi)),\ \xi\in\calo,$$where $\vf$ and $\vf_\vp$ are in \eqref{ee9.4}. We have
$$\na\phi_\vp=g'(\vf_\vp)\na\vf_\vp,\ \D\phi_\vp=g'(\vf_\vp)\D\vf_\vp+g''(\vf_\vp)|\na\vf_\vp|^2,\ \xi\in\calo,$$
and so, by \eqref{ee9.5},
\begin{equation}\label{e8.19}
\barr{lcl}
\dd\frac{\phi^2}{\phi_\vp}-\D\phi_\vp&=&
\dd g'(\vf_\vp)\(\frac{\vf^2}{\vf_\vp}-\D\vf_\vp\)+
\frac{g^2(\vf)}{g(\vf_\vp)}\vsp
&&
-g'(\vf_\vp)\,\dd\frac{\vf^2}{\vf_\vp}-g''(\vf_\vp)|\na \vf_\vp|^2\vsp
&\le& g'(\vf_\vp)|\na y|+\dd\frac{g^2(\vf)}{g(\vf_\vp)}-g'(\vf_\vp)\,\dd\frac{\vf^2}{\vf_\vp}.\earr\end{equation}

Now, proceeding as in the proof of Proposition \ref{p9.1}, we take $0\in\pp\calo$ and represent locally $\pp\calo$ as $\{(\xi',\xi_N);\ \xi_N=\g(\xi')\}$, where $\g\in C^2$, $\g(0)=0$, $\na\g(0)=0$. By \eqref{813} and since $g$ is increasing,  we have
$$D_N\phi_\vp(0)=
g'(\vf_\vp(0))D_N\vf_\vp(0)=-g'(\vf_\vp(0))
\,\dd\frac{\pp\vf_\vp}{\pp n}\,(0)\ge 0 .$$This yields
$\frac{\pp\phi_\vp}{\pp n}\,(0)=-D_N\phi_\vp(0)\le0$ and, therefore, replacing $0$ by an arbitrary point of $\pp\calo$,   we obtain that
$$\frac{\pp\phi_\vp}{\pp n}\le0\ \ \mbox{on }\pp\calo.$$Integrating \eqref{e8.19} over $\calo$, we  therefore get
$$\int_\calo\frac{\phi^2}{\phi_\vp}\,d\xi
\le\int_\calo\(g'(\vf_\vp)\(|\na y|-\frac{\vf^2}{\vf_\vp}\)+\frac{g^2(\vf)}{g(\vf_\vp)}\)d\xi.$$
Letting $\vp\to0$, we see that
$$\barr{l}
\dd\int_\calo g(\vf)d\xi\le\int_\calo g'(\vf)(|\na y|-\vf)d\xi+\dd\int_\calo g(\vf)d\xi\le\dd\int_\calo g(|\na y|)d\xi\earr$$
because $g'(u)(u-v)\ge g(u)-g(v),$ $\ff u,v\in\rr^+$. This completes the proof of \eqref{e8.18}. $\Box$\bk

Let $j_\lbb:\rr^N\to\rr$ be the Moreau--Yosida approximation from Section 5. Then, since $\na j_\lbb=\psi_\lbb$, it is easy to check that
$$j_\lbb(u)=\left\{\barr{ll}
\dd\frac1{2\lbb}\,|u|^2_N&\mbox{ for }|u|_N\le\lbb,\vsp
|u|_N-\dd\frac\lbb2&\mbox{ for }|u|_N>\lbb.\earr\right.$$

\begin{corollary}\label{c8.8} For all $\vp>0$ and $\lbb>0$, we have
\begin{equation}\label{e9.17}\int_\calo j_\lbb(\na J_\vp(y))d\xi\le\int_\calo j_\lbb(\na y(\xi))d\xi,\ \ff y\in H^1_0(\calo).\end{equation}
\end{corollary}

\n{\bf Proof.} One applies Proposition \ref{p8.6} to the function  $$\mbox{$g(r)=\left\{\barr{ll}
 \dd\frac1{2\lbb}\,r^2&\mbox{for }0\le r\le\lbb,\vsp
  r-\dd\frac\lbb2&\mbox{for }r>\lbb.\ \ \ \ \  \Box\earr\right.$}$$

\begin{remark}\label{r8.9} {\rm In Corollary \ref{c8.8}, the quadratic growth condition on $g$ can be relaxed. If, e.g., $g$ grows at most of order $p\in[1,\9)$, then
  $$\int g(|\na J_\vp(y)|)d\xi\le\int g(|\na y|)d\xi,\ \ff y\in W^{1,p}_0(\calo).$$ In particular, applying Corollary \ref{c8.8} to $g(u)=|u|^p$, where $1\le p<\9$, we~obtain that for each bounded and convex set $\calo\subset\rr^N$ with $C^2$-boundary, we have
\begin{equation}\label{e8.20}
|\na J_\vp(y)|_p\le|\na y|_p,\ \ff y\in W^{1,p}_0(\calo).\end{equation}
The case $p=\infty$ is also true and was earlier proved by  Brezis  and Stampacchia (\cite{15a}).
In other words, the operator $A$ is dissipative in $W^{1,p}_0(\calo)$ for all $1\leq p \leq \infty$.}\end{remark}

\section{Appendix 2. (Proof of \eqref{e511primm})}
\setcounter{equation}{0}

We have, for $y,z\in H^1_0(\calo)$,
 $$\barr{l}
 \<\wt A_\lbb(t)y-\wt A_\lbb(t)z,y-z\>\vsp
  =\dd\int_\calo(\psi_\lbb
 (\na(e^{W(t)}y))-\psi_\lbb(\na(e^{W(t)}z)))\cdot
 \na(e^{-W(t)}(y-z))d\xi\vsp
 +\lbb\dd\int_\calo\na(e^{-W(t)}\cdot(y-z))\cdot
 \na(e^{W(t)}(y-z))d\xi
 +\dd\frac12\int_\calo\mu|y-z|^2d\xi\vsp
 =\dd\int_\calo(\psi_\lbb(\na(e^{W(t)}y))-\psi_\lbb(\na(e^{W(t)}z)))
 \cdot(\na(e^{W(t)}y)-\na(e^{W(t)}z))
 e^{-2W(t)}d\xi\vsp
 +\dd\int_\calo(\psi_\lbb(\na(e^{W(t)}y))-\psi_\lbb(\na(e^{W(t)}z)))
 e^{W(t)}(y-z)\cdot\na(e^{-2W(t)})d\xi\vsp
 +\dd\frac12\int_\calo\mu|y-z|^2d\xi+\lbb\int_\calo\na(e^{-W(t)}(y-z))\cdot
 \na(e^{W(t)}(y-z))d\xi\vsp
 \ge-2{\rm Lip}\psi_\lbb
\dd \int_\calo|\na(e^{W(t)}(y-z))||y-z|\,|\na W(t)|\,e^{-W(t)}d\xi\vsp
+\lbb\dd\int_\calo|\na(y-z)|^2d\xi
-\lbb\int_\calo(y-z)^2|\na W(t)|^2d\xi\vsp
\ge-2{\rm Lip}\,\psi_\lbb\dd\int_\calo
(|\na(y-z)||y-z||\na W(t)|
+|y-z|^2|\na W(t)|^2)d\xi\vsp
+\lbb\dd\int_\calo(|\na(y-z)|^2
-(y-z)^2|\na W(t)|^2)d\xi
\ge-\de^\lbb_t(\oo)|y-z|^2_2,\ \ff t\in[0,T],\earr$$
where $$\de^\lbb_t(\oo)=\(\frac1\lbb({\rm Lip}\,\psi_\lbb)^2+2\)|\na W(t)(\oo)|^2_\9,$$ ${\rm Lip}\,\psi_\lbb$ is the Lipschitz constant of $\psi_\lbb$ and we have used the Young inequality in the last step. Then, \eqref{e511primm}  follows.

\end{document}